\documentclass[a4paper, 12pt]{article}
\usepackage{fullpage}
\usepackage{multirow}
\usepackage{amsmath, amsthm, amssymb, mathrsfs, graphicx, subfigure}
\usepackage{ifthen, url}
\usepackage[usenames]{color}
\usepackage{times}
\usepackage{bm}

\usepackage{natbib}

\theoremstyle{plain}

\theoremstyle{definition}

\newsavebox{\tablebox}
\newcommand{\prob}{{\Bbb P}}
\newcommand{\expn}{{\Bbb E}}

\newtheorem{thm}{Theorem}

\newtheorem{lem}[thm]{Lemma}
\newtheorem{remark}{Remark}
\newtheorem{example}{Example}

\begin{document}

\title{Optimal hybrid block bootstrap for sample quantiles under weak dependence}

\author{
Todd A. Kuffner \\
Department of Mathematics \\
Washington University in St. Louis \\
\url{kuffner@wustl.edu} \\
\mbox{} \\
Stephen M.S. Lee \\
Department of Statistics and Actuarial Science \\
The University of Hong Kong \\
\url{smslee@hku.hk} \\
\mbox{} \\
G. Alastair Young \\
Department of Mathematics \\
Imperial College London \\
\url{alastair.young@imperial.ac.uk}
}

\maketitle

\begin{abstract}
We establish a general theory of optimality for block bootstrap distribution estimation for sample quantiles under a mild strong mixing assumption. In contrast to existing results, we study the block bootstrap for varying numbers of blocks. This corresponds to a hybrid between the subsampling bootstrap and the moving block bootstrap (MBB), in which the number of blocks is somewhere between 1 and the ratio of sample size to block length. Our main theorem determines the optimal choice of the number of blocks and block length to achieve the best possible convergence rate for the block bootstrap distribution estimator for sample quantiles. As part of our analysis, we also prove an important lemma which gives the convergence rate of the block bootstrap distribution estimator, with implications even for the smooth function model. We propose an intuitive procedure for empirical selection of the optimal number and length of blocks. Relevant examples are presented which illustrate the benefits of optimally choosing the number of blocks.

\medskip

\end{abstract}

\section{Introduction}
\label{sec:intro}

A broad categorization of settings for bootstrap methods might be as follows: (i) smooth functionals of independent data; (ii) nonsmooth functionals of independent data; (iii) smooth functionals of dependent data; and (iv) nonsmooth functionals of dependent data. Setting (i) is the classic setting of the bootstrap, with \citet{Hall:1992} being an authoritative reference. In setting (ii), bootstrap methods for approximating distributions of sample quantiles have been studied by \citet{Efron:1979, BickelFreedman:1981, Singh:1981, Babu:1986, Efron:1982, GPSB:1984, HallSheather:1988, HDR:1989, HallMartin:1991, DEANGELIS:1993} and \citet{FalkJanas:1992}. 

In setting (iii), the existing literature is concentrated on block bootstrap methods for smooth functionals with dependent data, beginning with \citet{Hall:1985}, and \citet{Carlstein:1986}. Subsequently, the moving block bootstrap (MBB) was proposed by \citet{Kunsch:1989} and \citet{LiuSingh:1992}. Other variants and properties of block bootstrap methods were considered by \citet{PaparoditisPolitis:2001, PolitisRomano:1992, PolitisRomano:1994a} and \citet{PRW:1997}, to name a few. In setting (iii), various block bootstrap methods and their properties for weakly dependent sequences have been investigated by, for example, \citet{Buhlmann:1994, NNR:1994, HHJ:1995, GotzeKunsch:1996, Lahiri:1992, Lahiri:1996, Lahiri:1999} and \citet{BuhlmannKunsch:1999}. 

In contrast, there is a less extensive literature about block bootstrap methods in setting (iv), i.e. for nonsmooth functions of dependent data. \citet{SunLahiri:2006}, \citet{Sun:2007} and \citet{SharipovWendler:2013} are notable exceptions. Those authors considered block bootstrap approximation for sample quantiles under weak dependence. \citet{SunLahiri:2006} established strong consistency of the MBB, assuming only a polynomial (strong) mixing rate, for both distribution and variance estimation of the sample quantiles. \citet{SharipovWendler:2013} established similar results for the circular block bootstrap utilizing a different set of conditions to take advantage of empirical process theory for the Bahadur-Ghosh representation of the sample quantile. \citet{Sun:2007} is particularly relevant to our work, as discussed further below. All of these earlier results assume that the number of blocks tends to infinity with the sample size. 

More recently, \citet{KLY:2017a} established a more general consistency result for a hybrid block bootstrap, for both distribution and variance estimation of sample quantiles. While an exponential mixing rate is assumed, \citet{KLY:2017a} proved weak consistency for \textit{any} number of blocks, $1\le b=O(n/\ell)$ (as $n \rightarrow \infty$), whereas the existing proofs for the MBB and circular block bootstrap required that $b \rightarrow \infty$, where $b=\lfloor n/\ell \rfloor$. Here, $n$ is the available sample size, and $\ell$ is the block length. The value of $b$ is the number of resampled blocks to be pasted to form the bootstrap data series. The case $b=1$ corresponds to the subsampling bootstrap \citep{PolitisRomano:1994a}, and the case $b=\lfloor n/\ell \rfloor$ is the standard MBB \citep{Kunsch:1989}. Therefore, the consistency results in \citet{KLY:2017a} are for a `hybrid' between the MBB and the subsampling bootstrap, and those two extremes are covered by the same theory.

As noted in \citet{KLY:2017a}, their theoretical and empirical results suggest that there can be substantial performance improvement, in terms of mean squared errors (MSE) for both the variance and distribution estimators, when choosing some value of $b >1$, but less than $\lfloor n/\ell \rfloor$. This suggests the following question: does there exist some optimal choice of the pair $(b, \ell)$ which provides the best convergence rate for the bootstrap distribution estimator for sample quantiles under weak dependence? We answer that question in the present paper.

Related to the motivation of the present paper is the paper by \citet{Sun:2007}. She studied the convergence rate of the moving block bootstrap distribution estimator for sample quantiles with dependent data. A strong mixing condition with exponentially decaying mixing coefficients was assumed. An almost sure convergence result was established, and the best rate of convergence was found to be $O(n^{-1/4}\log \log n)$. We consider a weaker polynomial rate condition, which is also slightly weaker than that assumed in \citet{SunLahiri:2006}. Moreover, we allow the number of blocks to vary, instead of fixing $b= \lfloor n/\ell \rfloor$. Our main theorem establishes the convergence rate of the `hybrid' bootstrap distribution estimator for sample quantiles. It is a hybrid between the MBB ($b=\lfloor n/\ell \rfloor$) and subsampling ($b=1$) bootstrap. We also apply our theory to the setting of \citet{Sun:2007} below.

Another recent development in dependent data bootstrap methodology is the convolved subsampling bootstrap \citep{TPN:2017}. This bootstrap estimator is defined by the $k$-fold self-convolution of a subsampling distribution. In the special case of the sample means problem, this corresponds to our hybrid bootstrap. For the sample quantile problem which is the particular focus here, convolved subsampling bootstrap essentially computes the average of within-block sample quantiles over the $b$ resampled blocks. By contrast, our estimator is the sample quantile of a single series formed by joining $b$ blocks. Further theoretical comparison of these approaches will be undertaken elsewhere. 

Other indirectly-related work includes \citet{Lahiri:2005}, who studied consistency of jackknife-after-bootstrap variance estimation for bootstrap quantiles, and \citet{GLN:2015}, who showed that the \citet{SunLahiri:2006} strong consistency results for distribution and variance estimation via the MBB also hold for the smoothed extended tapered block bootstrap (SETBB). We mention that \citet{ShaoPolitis:2013} employed a fixed $b$ subsampling procedure to estimate confidence sets for statistics adhering to the smooth function model, which of course does not include sample quantiles.

Aside from our general optimality results being of foundational and practical value, they also indicate that adaptive selection of the number of blocks could yield considerable improvements in convergence rates for block bootstrap distribution estimators. Moreover, Lemma~\ref{lem:bootcdf} below is of independent interest, as it gives the convergence rate of the block bootstrap distribution estimator, and has bearing on the regular smooth function model. We have included several relevant empirical examples to illustrate the potential gains of optimal choice of the number of blocks, as opposed to using the prescribed value of $b$ for either the subsampling bootstrap $(b=1)$ or the MBB ($b=\lfloor n/\ell \rfloor$). In \S~\ref{sec:empirical}, we give practical guidance as to how to choose $(b,\ell)$ in a given applied problem, by proposing a procedure for this purpose.

\section{Problem Setting}
\label{sec:setting}

Let $\mathbb{Z} \equiv \{0, \pm 1, \pm 2, \ldots \}$ be the set of all integers. Define $\{X_{i}\}_{i \in \mathbb{Z}}$ to be a doubly-infinite sequence of random variables on the probability space $(\Omega, \mathcal{F}, P)$. The elements of the sequence possess a common distribution function $F$, and its corresponding quantile function $F^{-1}$, defined by 
\[
F^{-1}(p)=\inf\{u : F(u) \geq p\}, \qquad p \in (0,1).
\]
We will study the block bootstrap distribution estimator of a suitably centered and scaled sample quantile. It is assumed throughout that $\{X_{i}\}_{i \in \mathbb{Z}}$ is a strictly stationary process. The sequence $(X_{1}, \ldots, X_{n})$ denotes a sample of size $n$ from $\{X_{i}\}_{i \in \mathbb{Z}}$. 

\subsection{The Block Bootstrap}

The moving blocks bootstrap (MBB) \citep{Kunsch:1989} splits the original sample $(X_{1}, \ldots, X_{n})$ into overlapping blocks of size $\ell$, $B_{i}=(X_{i}, \ldots, X_{i+\ell-1})$, together constituting a set \break$\{B_{1}, \ldots, B_{n-\ell+1}\}$. Let $B_{1}^{*}, \ldots, B_{b}^{*}$ be a random sample drawn with replacement from the original blocks, where $b=\lfloor n/\ell \rfloor$ is the number of blocks that will be pasted together to form a pseudo-time series. For a real number $h$, the notation $\lfloor h \rfloor $ is defined as the largest integer $\leq h$, and $\lceil h \rceil$ is the smallest integer $\geq h$. That $B_{1}^{*}, \ldots, B_{b}^{*}$ is a random sample from $\{B_{1}, \ldots, B_{n-\ell+1}\}$ means that the sampled blocks are independently and identically distributed according to a discrete uniform distribution on $\{B_{1}, \ldots, B_{n-\ell+1}\}$. The observations in the $i$th resampled block, $B_{i}^{*}$, are $X_{(i-1)\ell+1}^{*}, \ldots, X_{i\ell}^{*}$, for $1 \leq i \leq b$. Then the MBB sample is the concatenation of the resampled blocks, written as 
\[
\underbrace{X_{1}^{*}, \ldots, X_{\ell}^{*}}_{B_{1}^{*}}, \underbrace{X_{\ell+1}^{*}, \ldots, X_{2\ell}^{*}}_{B_{2}^{*}}, \underbrace{X_{2\ell+1}^{*}, \ldots, X_{(b-1)\ell}^{*}}_{B_{3}^{*} \ldots B_{(b-1)}^{*}}, \underbrace{X_{(b-1)\ell+1}^{*}, \ldots, X_{b\ell}^{*}}_{B_{b}^{*}}.
\]
Note that this way of constructing the pseudo-time series will reproduce the original dependence structure \textit{asymptotically}.

The subsampling bootstrap \citep{PolitisRomano:1994a}, and specifically the overlapping blocks version relevant to the present setting, first splits the original sample into precisely the same overlapping blocks as the MBB, each of length $\ell$. However, the subsampling bootstrap draws only a single block. A nice property of this procedure is that the original dependence structure in the sample is exactly retained in the single subsample. By contrast, the pseudo-time series constructed by the MBB only reproduces the original dependence structure asymptotically. 

We define dependence for the sequence of random variables $\{X_{i}\}_{i \in \mathbb{Z}}$ in terms of the mixing properties of $\sigma$-algebras generated by subsets of the sequence which are separated by a distance, in units of time, tending to infinity. For any two sub-$\sigma$-algebras of $\mathcal{F}$, say $\mathcal{F}_{1}$ and $\mathcal{F}_{2}$, the $\alpha$-mixing coefficient between $\mathcal{F}_{1}$ and $\mathcal{F}_{2}$ is defined to be \citep[Section 16.2.1]{AthreyaLahiri:2006}
\begin{equation}
\label{eqn:alpha}
\alpha(\mathcal{F}_{1}, \mathcal{F}_{2}) \equiv \sup_{A \in \mathcal{F}_{1}, B \in \mathcal{F}_{2}} \vert P(A \cap B)-P(A)P(B) \vert .
\end{equation}
Write $\mathcal{F}_{k}^{k+t}$ for the smallest $\sigma$-algebra of subsets of $\Omega$ with respect to which $X_{i}$, $i=k, \ldots, k+t$, are measurable. Let $\mathcal{F}_{-\infty}^{k}$ be the smallest $\sigma$-algebra which contains the unions of all of the $\sigma$-algebras $\mathcal{F}_{a}^{k}$ as $a \rightarrow -\infty$. That is, $\mathcal{F}_{-\infty}^{k}$ is a sub-$\sigma$-algebra of $\mathcal{F}$, and it is the $\sigma$-algebra generated by the random variables $X_{a}, X_{a+1}, \ldots, X_{k}$ as $a\rightarrow -\infty$. Similarly, for $-\infty \leq k  \leq \infty$, let $\mathcal{F}_{k}^{\infty}$ be the $\sigma$-algebra generated by the random variables $X_{k+1}, X_{k+2}, \ldots, X_{k+a}$, as $a \rightarrow \infty$. The $\alpha$-mixing coefficient of the sequence $\{X_{i}\}_{i \in \mathbb{Z}}$ is defined as
\[
\alpha(t) \equiv \sup_{k \in \mathbb{Z}} \alpha(\mathcal{F}_{-\infty}^{k}, \mathcal{F}_{k+t}^{\infty}),
\]
where $\alpha(\cdot, \cdot)$ is defined in \eqref{eqn:alpha}. If the $\alpha$-mixing coefficient decays to zero, 
\begin{equation}
\label{eqn:alphamix}
\lim_{t \rightarrow \infty} \alpha(t) = 0,
\end{equation}
then the process $\{X_{i}\}_{i \in \mathbb{Z}}$ is said to be strongly mixing. The sequence of random variables $\{X_{i}\}_{i \in \mathbb{Z}}$ is said to be weakly dependent if the process $\{X_{i} \}_{i \in \mathbb{Z}}$ is strongly mixing, i.e. if \eqref{eqn:alphamix} holds.

\section{Theoretical Results}
\label{sec:theory}

Assume that $(X_1,\ldots,X_n)$ is a sample of a stationary
strong mixing process with mixing
coefficient $\alpha(t)$ satisfying, for some $\beta>5$,
\[ \alpha(t)=O(t^{-\beta}),\qquad t\rightarrow\infty. \]
Denote by $F$ the distribution function of $X_1$ and
$F_n$ the empirical distribution function of $(X_1,\ldots,X_n)$.

Define, for $\ell\in\{1,2,\ldots,n\}$, $b\in\{1,2,\ldots\}$ and $x\in\mathbb{R}$,
$J_1,\ldots,J_b$ to be independent random indices uniformly drawn from the
set $\{1,\ldots,n-\ell+1\}$,
\begin{align*}
U_i(x)&=\ell^{-1}\sum_{t=i}^{i+\ell-1}\pmb{1}\{X_t\le x\},\;\;
i=1,\ldots,n-\ell+1,\\
U^*_i(x)&=\ell^{-1}\sum_{t=J_i}^{J_i+\ell-1}\pmb{1}\{X_{t}\le x\},\;\;
i=1,\ldots,b,
\end{align*}
\begin{equation*}
    \tilde{F}_n(x)=(n-\ell+1)^{-1}\sum_{i=1}^{n-\ell+1}U_i(x), \quad F^*_n(x)=b^{-1}\sum_{i=1}^bU^*_i(x).
\end{equation*}
Define, for $p\in(0,1)$,
\[ \xi_p=F^{-1}(p),\qquad\hat\xi_n=F_n^{-1}(p),\qquad\tilde\xi_n=\tilde{F}_n^{-1}(p),
\qquad\xi^*_n=F^{*-1}_n(p).\]
Assume that $f=F'$ is defined on a neighbourhood $\mathscr{N}_p$ of $\xi_p$,  with
\[0<\underset{x\in\mathscr{N}_p}\inf f(x)\le \underset{x\in\mathscr{N}_p}\sup f(x)<\infty.\]
\begin{thm}
\label{thm:main}
Suppose that $n^{-\frac{4\beta+7}{6(3\beta+5)}}\ell\rightarrow\infty$ and
$b\ge 1$.
Let $x\in\mathbb{R}$ be fixed and $\delta>0$ be any arbitrarily small constant.
\begin{itemize}
\item[(i)] If $\beta\in(5,\infty)$ and $\ell=O(b)$, then
\begin{eqnarray*}
\lefteqn{
\prob\Big(\sqrt{b\ell}\big(\xi^*_n-\tilde\xi_n\big)\le x
\Big|X_1,\ldots,X_n\Big)}\nonumber\\
&=&\prob\Big(\sqrt{n}\big(\hat\xi_n-\xi_p\big)\le x\Big)
+
O_p\Big(\ell^{-1}+\ell^{1/2}n^{-1/2}+(b\ell)^{-1/2}\ell^\delta\nonumber\\
&&\qquad +\,n^{-\frac{\beta-1}{2(\beta+1)}+\delta}(b\ell)^{(1-\delta)/4}
+n^{-1}b^{\frac{2\beta+1}{4(\beta+1)}-\delta}
\ell^{\frac{4\beta+7}{4(\beta+1)}+5\delta}
\Big)\nonumber\\
&&+\,
o_p\Big(
n^{-\frac{\beta-3}{\beta-1}+\delta}(b\ell)^{1/2}
+n^{-\frac{3\beta-1}{4(\beta+1)}+\delta}(b\ell)^{1/2} \nonumber\\
&&\qquad+\,
n^{-\frac{\beta(2\beta-3)}{(\beta-1)(2\beta+1)}+\delta}b^{\frac{1}{2}}
\ell^{\frac{1}{2}+\frac{2(\beta+3)}{(\beta-1)(2\beta+1)}}
+
n^{-\frac{4\beta+5}{4(\beta+1)}+\delta}b^{\frac{1}{2}}
\ell^{\frac{\beta+2}{\beta+1}}\nonumber\\
&&\qquad+\,
n^{-\frac{2(\beta+1)}{2\beta+1}+\delta}b^{\frac{1}{2}}
\ell^{\frac{4\beta+7}{2(2\beta+1)}}
+n^{-\frac{4\beta^2+3\beta+1}{2(2\beta+1)(\beta+1)}+\delta}
b^{\frac{1}{2}}\ell^{\frac{3\beta+4}{2(2\beta+1)}}\Big).
\end{eqnarray*}
\item[(ii)]
If $\beta=\infty$, then
\begin{eqnarray*}
\lefteqn{
\prob\Big(\sqrt{b\ell}\big(\xi^*_n-\tilde\xi_n\big)\le x
\Big|X_1,\ldots,X_n\Big)}\nonumber\\
&=&\prob\Big(\sqrt{n}\big(\hat\xi_n-\xi_p\big)\le x\Big)
+
O_p\Big(\ell^{-1}+\ell^{1/2}n^{-1/2}+(b\ell)^{-1/2}\nonumber\\
&&\qquad+\,
n^{-1}b^{\frac{1}{2}-\delta}
\ell^{1+5\delta}
+n^{-\frac{1}{2}+\delta}(b\ell)^{(1-\delta)/4}\Big)\nonumber\\
&&\qquad+\,
o_p\Big(
n^{-\frac{3}{4}+\delta}(b\ell)^{1/2}+n^{-1+\delta}b^{\frac{1}{2}}\ell\Big).
\end{eqnarray*}
\end{itemize}
\end{thm}
We may deduce from Theorem~\ref{thm:main} the following:
\begin{description}
\item[Case (i).]  $\beta\in(5,\infty)$ and $\ell=O(b)$.

The convergence rate of the bootstrap distribution estimator is minimised
by setting
\[\ell\propto b\propto
\begin{cases}
n^{\frac{4\beta+7}{6(3\beta+5)}}\log n,&\beta\in\big(5,(7+\sqrt{185})/4\big],\\
n^{\frac{\beta-1}{3(\beta+1)}},&\beta\in\big(7+\sqrt{185})/4,\infty\big),
\end{cases}
\]
which yields, for any $\delta>0$,
\begin{eqnarray}
\lefteqn{
\prob\Big(\sqrt{b\ell}\big(\xi^*_n-\tilde\xi_n\big)\le x
\Big|X_1,\ldots,X_n\Big)-
\prob\Big(\sqrt{n}\big(\hat\xi_n-\xi_p\big)\le x\Big)}\nonumber\\
&=&
\begin{cases}
O_p\left(n^{-\frac{14\beta^2+\beta-37}{12(3\beta+5)(\beta+1)}+\delta}\right),
&\beta\in\big(5,(7+\sqrt{185})/4\big],\\
O_p\left(n^{-\frac{\beta-1}{3(\beta+1)}+\delta}\right),
&\beta\in\big(7+\sqrt{185})/4,\infty\big).
\end{cases}
\label{oracle.rate}
\end{eqnarray}
Note that as $\beta\rightarrow\infty$, the optimal orders of $\ell$ and $b$ approach $n^{1/3}$, 
which does not depend on unknown parameters and
may be taken as a practical reference for empirical choices of $\ell$ and $b$.
With such choices, that is $\ell\propto b\propto n^{1/3}$, the bootstrap distribution estimator
has the convergence rate
$O_p\left(n^{-\frac{\beta-2}{3(\beta+1)}+\delta}\right)$, for $\beta\in(5,\infty)$
and any $\delta>0$. The latter convergence rate is slightly slower than that specified in (\ref{oracle.rate}), a price to pay for
the absence of knowledge of $\beta$.

On the other hand, the MBB sets $b=\lfloor n/\ell\rfloor$, based on which the optimal $\ell$ is of order $n^{1/3}$, so that $b\propto n^{2/3}$.
The convergence rate of the resulting bootstrap distribution estimator is given, for any $\delta>0$, by
\[
\begin{cases}
O_p\left(n^{-\frac{\beta-5}{2(\beta-1)}+\delta}\right),
&\beta\in\big(5,2+\sqrt{17}\,\big],\\
O_p\left(n^{-\frac{\beta-3}{4(\beta+1)}+\delta}\right),
&\beta\in\big(2+\sqrt{17},\infty\big),
\end{cases}
\]
which is markedly slower than that obtained by setting $\ell\propto b\propto n^{1/3}$.

Figure~\ref{errorplot} compares the optimal convergence rate with those based on $b\propto\ell\propto n^{1/3}$ and
$b=\lfloor n/\ell\rfloor\propto n^{2/3}$, respectively.
\begin{figure}[htp]
\centering
\caption{\label{errorplot}
From Theorem~\ref{thm:main}, Case (i). Log error rates for the block bootstrap distribution estimator are plotted against $\beta$ for the optimal pairs of $(b,\ell)$. The choice $b=\ell=n^{1/3}$ is optimal for $\beta=\infty$, and it is our recommendation when no information about the exact value of $\beta$ is available. Thus the discrepancy between the blue and red curves shows how ignorance about $\beta$ affects the error rate. The MBB choice $(b=n^{2/3},\ell=n^{1/3})$ is optimal for K\"{u}nsch's MBB.}
\includegraphics[scale=0.7]{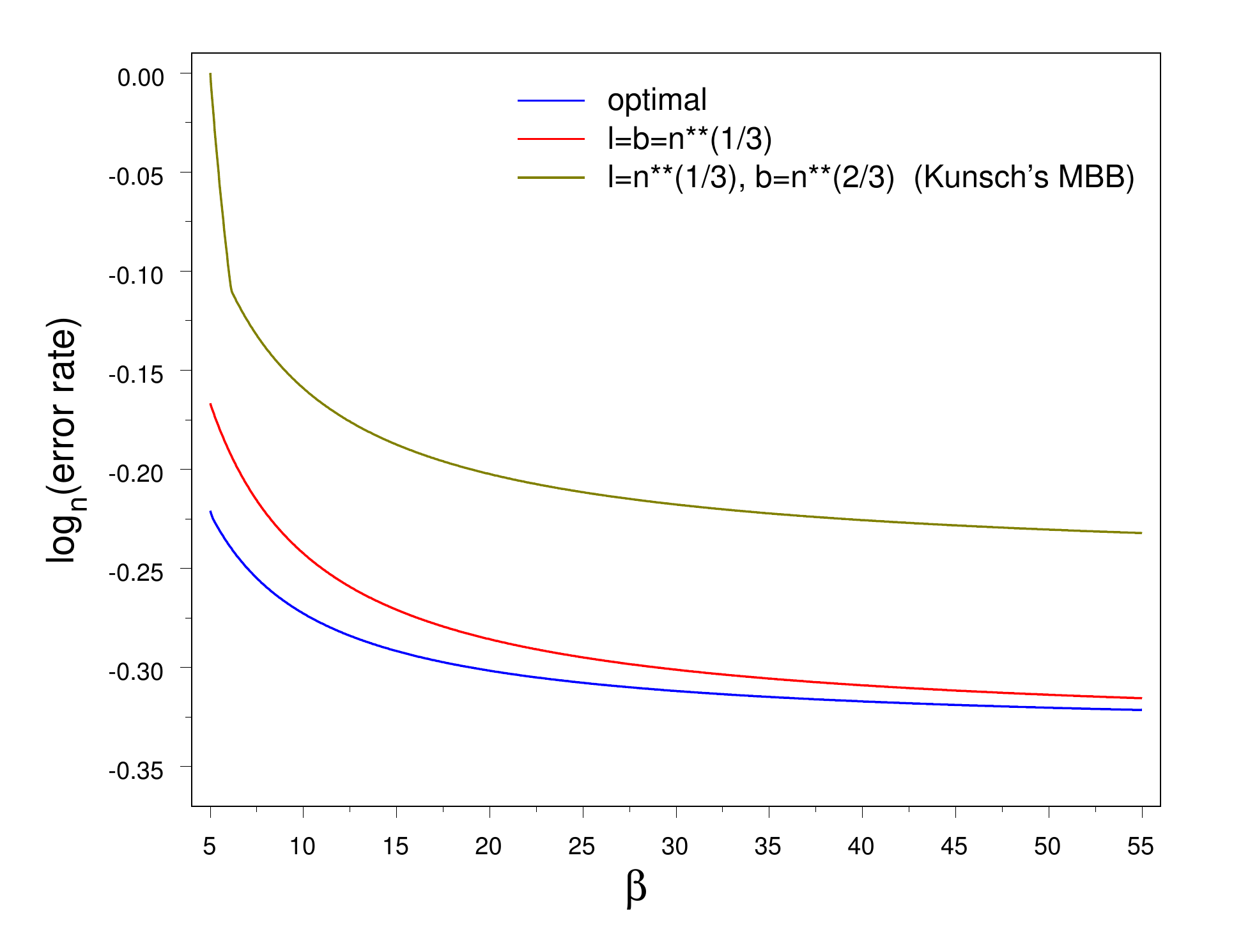}
\end{figure}
\item[Case (ii).] $\beta=\infty$.

The error rate has an order minimised by setting $\ell\propto b\propto n^{1/3}$, which yields
\[
\prob\Big(\sqrt{b\ell}\big(\xi^*_n-\tilde\xi_n\big)\le x
\Big|X_1,\ldots,X_n\Big)=\prob\Big(\sqrt{n}\big(\hat\xi_n-\xi_p\big)\le x\Big)
+O_p\Big(n^{-1/3+\delta}\Big),\]
for any arbitrarily small $\delta>0$.
For MBB, the error rate is minimised if $\ell$ is chosen to have order between $n^{1/4}$ and $n^{1/2}$, yielding an optimal convergence rate
of order $O_p\left(n^{-1/4+\delta}\right)$ for any $\delta>0$.
If we set $b=1$, which amounts to the subsampling method, then the fastest error rate has order $O_p\big(n^{-1/4}\big)$,
attained by setting $\ell\propto n^{1/2}$.

\end{description}

\begin{remark}
\label{remark:Sun}
Case (ii) applies to any strong mixing setting where the mixing coefficients decay faster than any finite polynomial rate. Namely, this includes the exponential rate of \citet{Sun:2007}, and anything faster. The case that $\beta=\infty$ could be replaced by an assumption that the mixing coefficient satisfies $\alpha(t)=O(e^{-Ct})$ as $t\rightarrow \infty$, for some $C>0$. This is the assumption in \cite{KLY:2017a}, which is equivalent to the condition in \citet{Sun:2007}. Under such an exponential rate condition, the results for Case (ii) would still hold. However, this would unnecessarily exclude any rate that is faster than polynomial but slower than exponential. We emphasize that \citet{Sun:2007} proved an almost sure convergence result under this exponential rate condition, using only the MBB choices of $(b,\ell)$.
\end{remark}

\begin{remark}
\label{remark:lemma}
Of independent interest is the result of Lemma~\ref{lem:bootcdf}, which gives the convergence rate
 of the block bootstrap distribution estimator for $\sqrt{n}\big(F_n(x)-F(x)\big)$ and has a bearing on the regular smooth function model. 
 Consider the simpler case where $\beta=\infty$. It is easily seen that the convergence rate is minimised at $O_p(n^{-1/3})$,
 attained by setting 
 $\ell\propto n^{1/3}$ and $b$ having order not smaller than $n^{1/3}$, of which MBB is a special case.
 The subsampling method ($b=1$), however, has at best a convergence rate of only order $O_p(n^{-1/4})$, attained
 by setting $\ell\propto n^{1/2}$.
 \end{remark}

\begin{remark}
Results on distribution estimation for $\sqrt{n}\big(\hat\xi_n-\xi_p\big)$, embodied in Theorem~\ref{thm:main}, differ substantially
from the regular case in that local estimation of $F$ over a shrinking neighbourhood of size $O_p\big((b\ell)^{-1/2}\big)$ around
$\xi_p$ incurs an error of order $n^{-1/2}(b\ell)^{1/4}$, which favours a small $b$ and precludes MBB
from yielding an optimal convergence rate.
\end{remark}

\section{Relevance to Coverage Error}
\label{cov.error}
Define 
\begin{equation*}
    \hat{G}_{n}(x)=\prob\left(\sqrt{b\ell}(\xi_{n}^{*}-\tilde{\xi}_{n})\leq x \vert X_{1}, \ldots, X_{n} \right)
    \end{equation*}
    and let $\Delta(n,b,\ell)$ be defined by
\begin{equation} \label{delta-coverage}
\hat{G}_{n}(x)=\Phi \left(xf(\xi_{p})/\sigma(\xi_{p})\right) + \Delta(n, b, \ell) .
\end{equation}
Our main results in \S~\ref{sec:theory} establish the asymptotic order of $\Delta(n,b,\ell)$ and derive the optimal orders of $(b,\ell)$ which minimise that order.

A level $\alpha$ lower percentile confidence interval for $\xi_{p}$ is given by
\begin{equation*}
[\hat{\xi}_{n}-n^{-1/2}\hat{G}_{n}^{-1}(\alpha), \infty).
\end{equation*}
Noting from \eqref{delta-coverage} that
\begin{equation*}
\hat{G}_{n}^{-1}(\alpha)=\Phi^{-1}(\alpha)\sigma(\xi_{p})/f(\xi_{p})+O_{p}\left(\Delta(n, b, \ell)\right),
\end{equation*}
and using \eqref{pf:true} from \S~\ref{sec:proofs}, which follows from standard asymptotic properties of $F_{n}$ for dependent data \citep{LahiriSun:2009}, we obtain that
\begin{align*}
\prob\left(\xi_{p} \geq \hat{\xi}_{n}-n^{-1/2}\hat{G}_{n}^{-1}(\alpha)\right)&=\prob\left(\sqrt{n}(\hat{\xi}_{n}-\xi_{p})\leq \Phi^{-1}(\alpha)\sigma(\xi_{p})/f(\xi_{p})\right) \\
& \qquad +O\left(\Delta(n,b,\ell)\right) \\
&= \alpha + O\left(\Delta(n,b,\ell)+n^{-1/2}\right) \\
&= \alpha + O\left(\Delta(n,b,\ell)\right) .
\end{align*}
Thus, minimising the order of $\Delta(n,b,\ell)$ also minimises the order of the coverage error.

\section{Practical Procedure for Selecting Optimal $(b,\ell)$}
\label{sec:empirical}

Setting $b=\lfloor c_{1}n^{1/3}\rfloor$ and $\ell= \lfloor c_{2}n^{1/3} \rfloor$, the objective is to find the optimal pair of positive constants $(c_{1}, c_{2})$ which minimise the estimation error of $\hat{G}_{n}(x)$, or coverage error under some obvious modification of the procedure. Note from \eqref{delta-coverage} and \eqref{pf:true} that
\begin{equation} \label{empirical}
\hat{G}_{n}(x)-G_{n}(x)=\Delta \left(n, \lfloor c_{1}n^{1/3} \rfloor, \lfloor c_{2}n^{1/3} \rfloor \right)\{1+o_{p}(1)\} .
\end{equation}
Define, for $c_{1}, c_{2}>0$ and a fixed $\rho>0$,
\begin{equation*}
\delta_{n}(c_{1}, c_{2})=\expn\left| \Delta \left(n, \lfloor c_{1}n^{1/3} \rfloor,  \lfloor c_{2}n^{1/3} \rfloor \right)\right|^{\rho}.
\end{equation*}
Then the $L_{\rho}$ estimation error of $\hat{G}_{n}(x)$ has the expansion
\begin{equation} \label{Lrho-error}
\expn\left|\hat{G}_{n}(x)-G_{n}(x)\right|^{\rho}=\delta_{n}(c_{1},c_{2})\{1+o(1)\} .
\end{equation}
We wish to minimise $\delta_{n}(c_{1},c_{2})$ with respect to $c_{1}, c_{2}$.

Let $M$ be a subsample size satisfying $M=o(n)$ and $M \rightarrow \infty$. Let $\hat{G}_{M}^{(j)}(x)$ be constructed analogously to $\hat{G}_{n}(x)$, with the complete sample $(X_{1}, \ldots, X_{n})$ replaced by the $j$th block of $M$ consecutive observations drawn from $(X_{1}, \ldots, X_{n})$, for $j=1, n-M+1$. Then we have, analogous to \eqref{empirical}, that
\begin{equation} \label{empirical-error}
\hat{G}_{M}^{(j)}(x)-G_{M}(x)=\Delta^{(j)}\left(M, \lfloor c_{1}M^{1/3} \rfloor, \lfloor c_{2}M^{1/3} \rfloor \right)\{1+o_{p}(1)\},
\end{equation}
where $\Delta^{(j)}(\cdot)$ denotes the version of $\Delta(\cdot)$ obtained from the $j$th subsample. Define
\begin{equation*}
Err(c_{1},c_{2})=(n-M+1)^{-1}\sum_{j}\left|\hat{G}_{M}^{(j)}(x)-\hat{G}_{n}(x)\right|^{\rho}.
\end{equation*}
Using \eqref{pf:true}, \eqref{empirical} and \eqref{empirical-error}, we have
\begin{align*}
\hat{G}_{M}^{(j)}(x)-\hat{G}_{n}(x) &= G_{M}(x)+\Delta^{(j)}\left(M,\lfloor c_{1}M^{1/3} \rfloor , \lfloor c_{2}M^{1/3} \rfloor \right)\{1+o_{p}(1)\} \\
& \qquad - G_{n}(x)-\Delta\left(n, \lfloor c_{1}n^{1/3} \rfloor , \lfloor c_{2}n^{1/3} \rfloor \right)\{1+o_{p}(1)\} \\
&= \Delta^{(j)}\left(M, \lfloor c_{1}M^{1/3} \rfloor , \lfloor c_{2}M^{1/3} \rfloor \right)\{1+o_{p}(1)\} .
\end{align*}
It follows that
\begin{eqnarray*}
Err(c_{1},c_{2})&=&(n-M+1)^{-1}\sum_{j}\left|\Delta^{(j)}\left(M, \lfloor c_{1}M^{1/3} \rfloor , \lfloor c_{2}M^{1/3} \rfloor \right)\right|^{\rho}\{1+o_p(1)\} \\&=&\delta_{M}(c_{1},c_{2})\{1+o_p(1)\}.    
\end{eqnarray*}
If we assume, as is typical, that $\delta_{n}(c_{1},c_{2})$ has a leading term of the form $\beta(c_{1},c_{2})n^{-\gamma}$ for some function $\beta(\cdot)$ independent of $n$, then $Err(c_{1},c_{2})$, $\delta_{M}(c_{1},c_{2})$ and $\delta_{n}(c_{1},c_{2})$ are all minimised at asymptotically the same $(c_{1},c_{2})$. Thus, an empirical procedure for choosing $(c_{1},c_{2})$, and hence choosing $(b, \ell)$, may be based on the minimisation of $Err(c_{1},c_{2})$.

This procedure constructs the error estimate $Err(c_1,c_2)$ by considering all $n-M+1$ subsamples of $M$ consecutive points drawn from the original data sample, and is therefore computationally expensive. However, the argument supporting minimization of this quantity actually only requires that the number of subsamples used in the construction should grow with sample size $n$. In practice, therefore, it is reasonable to evaluate the error measure $Err(c_1,c_2)$ using a smaller set of subsamples: in the numerical illustration given below, 20 subsamples, equally spaced along the data series $(X_1, \ldots, X_n)$, are used, allowing rapid evaluation of the error estimate.

\section{Examples}
To illustrate the benefits of optimally choosing $(b,\ell)$, we consider three very general examples. For concreteness, we consider $p=1/2$, and simulate the mean squared errors (MSEs) of hybrid block bootstrap estimators of $G_{n}(u)$ for particular choices of $u$. The true reference values of $G_{n}(\cdot)$ are approximated via massive simulation ($5 \times 10^{6}$ replications). For each of the sample sizes $n=200$, $n=500$, and $n=1000$, all entries in the included tables and heatmaps are based on $20,000$ replications, with $20,000$ bootstrap samples used within each replication, unless otherwise stated. For $n=2,000$, the number of replications and bootstrap samples are each $10,000$. For convenience, Table~\ref{values} provides some reference values of $(b,\ell)$ for MBB for the sample sizes we consider. This facilitates comparison with the MBB choice of $b=\lfloor n/\ell \rfloor$ for a range of values of $\ell$. In particular, we give values for $\ell$ approximately equal to $n^{1/2}$ (not optimal), $n^{1/3}$ (thought to be optimal), $n^{1/4}$, and $n^{1/5}$. 

\begin{table}[h]
\caption{\label{values}
Standard choices of $(b,\ell)$ for different $n$, with the MBB choice $b=\lfloor n/\ell \rfloor$.}
\begin{tabular}{l | c c c c}
& $(b, \ell \approx n^{1/2}$ & $(b, \ell \approx n^{1/3})$ & $(b, \ell \approx n^{1/4})$ & $(b, \ell \approx n^{1/5})$ \\ \hline
$n=200$ & $(14, 14)$ & $(33 , 6)$ & $(50 , 4)$ & $( 66, 3)$ \\
$n=500$ & $(22, 22)$ & $(62, 8)$ & $(100, 5)$ & $(125, 4)$ \\
$n=1,000$ & $(31, 32)$ & $(100, 10)$ & $(166, 6)$ & $(250, 4)$ \\
$n=2,000$ & $(44, 45)$ & $(153, 13)$ & $(285, 7)$ & $(400, 5)$
\end{tabular}
\end{table}

\begin{figure}[h]
\centering
\caption{\label{ARMAheat}
Heatmap for the ARMA(1,1) model with $n=200$.}
\includegraphics[height=9cm]{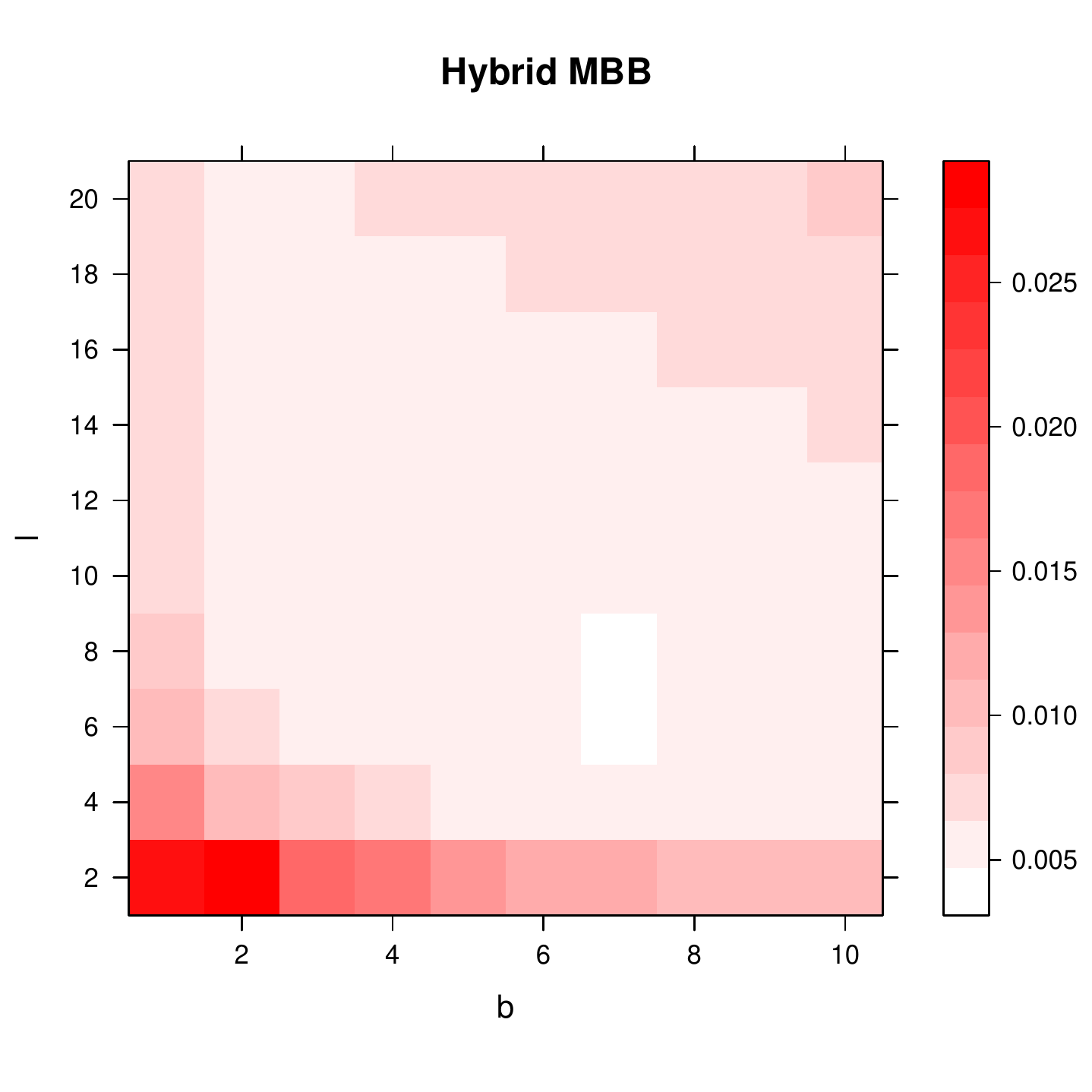}
\end{figure}

\begin{example}[ARMA(1,1)] Suppose that the observations are generated according to an ARMA~(1,1) model
\[
X_{t}-0.4 X_{(t-1)}=\epsilon_{t}+0.3 \epsilon_{(t-1)},
\]
with $\epsilon_{t}$ independent, identically distributed $N(0,1)$. The strong mixing condition is satisfied with an exponential rate \citep[Example 6.1]{Lahiri:2003}. An initial $X_{0}$ is sampled according to the marginal distribution, i.e. $X_{0} \sim N(0, 1.5833)$, and $\epsilon_{0} \sim N(0,1)$.
\end{example}
With $p=1/2$, we have $\xi_{p} =0$. We simulate the MSE in estimation of $G_n(1)$ over a range of $(b,\ell)$. The true value being estimated was computed (by massive simulation, as described) as $G_{n}(1) \approx 0.67978$. The heat map in Figure~\ref{ARMAheat} plots MSE for $n=200$, over a grid of values of $(b,\ell)$. The heat map clearly illustrates the sub-optimality of $b=1$, the subsampling bootstrap. The minimum MSE is 0.00472, with $(b,\ell)=(6,8)$. By contrast, the minimum MSE for the MBB is 0.00637, with $(b, \ell)=(33,6)$, and the subsampling bootstrap, which fixes $b=1$, has minimum MSE of 0.00754, with $\ell=14$.

We also compute the values of the pair $(b,\ell)$ which minimize MSE for other sample sizes, $n=500, 1000$, and $2000$. These results are shown in Table~\ref{ARMAtable}. Comparing with Table~\ref{values}, we note that the MSE-minimizing pair $(b,\ell)$ for each $n$ uses an $\ell$ strictly greater than $n^{1/3}$ and a $b$ much less than $\lfloor n/\ell \rfloor$. Additionally, the MSE-minimizing value of $b$ is much larger than $1$.

\begin{table}[h]
\caption{\label{ARMAtable}
ARMA(1,1) model. Choices of $(b,\ell)$ which minimize the MSE for estimating $G_{n}(1)$ for different sample sizes $n$.}
\centering
\begin{tabular}{l l l}
& $(b,\ell)$ & MSE \\ \hline
$n=200$  &   (6,8) &  0.00472 \\
$n=500$  &  (10,10) & 0.00250 \\
$n=1,000$ & (10,14) & 0.00154 \\ 
$n=2,000$ & (12,18) & 0.00097 \\
\end{tabular}
\end{table}

The theory says that the hybrid MBB has an error rate in estimation of $G_n(1)$ of $n^{-1/3}$, so we should expect the MSE to decrease at rate $n^{-2/3}$. In fact, a regression of $\log(MSE)$ on $\log(n)$ for the values reported in Table~\ref{ARMAtable} has slope $-0.6885$, which is not far off $-2/3$. The heatmap illustrates that the subsampling and MBB choices of $(b,\ell)$ are suboptimal.

For the current problem, of estimation of the sampling distribution of the sample quantile, there is therefore clear theoretical and practical advantage in using the hybrid block bootstrap, $b \ell <n, b \neq 1$, over the moving block bootstrap. Remark~\ref{remark:lemma} indicates, by contrast, that we might expect to see little
difference, in estimation error terms, between the hybrid block bootstrap procedure and MBB if, instead, we are interested in estimation of $\prob\{\sqrt{n}\big(F_n(x)-{F}(x)\big)\le y \}$. This was verified by considering, for all combinations of $(b, \ell)$, the MSE of the estimator $\prob\{\sqrt{b\ell}\big(F^*_n(x)-\tilde{F}_n(x)\big)\le y\Big|X_1,\ldots,X_n\}$, for $x=0$, so that $F(x)=0.5$, and $y=0.9$, for which the quantity being estimated $\approx 0.89501$, for sample size $n=100$. Based on 20,000 replications, with 20,000 bootstrap samples being used in construction of the estimator for each, the minimum MSE achieved by MBB is 0.00084, with $(b, \ell)=(25,4)$. This is very similar to the overall minimum MSE of 0.00082, seen for $(b, \ell)=(18,5)$.  The minimum MSE of the subsampling bootstrap, $b=1$, is 0.00334, substantially larger, when $\ell=7$. This same picture was seen for $n=200$, when, for the same values $x=0, y=0.9$, the true probability being estimated $\approx 0.87781$. Simulation shows that the minimum MSE of MBB is then 0.00108, with $(b, \ell)=(28,7)$, with the same minimum MSE for the hybrid block bootstrap, achieved for $(b, \ell)=(30,6)$. Here the subsampling bootstrap yields an optimal MSE of 0.00227 when $\ell=8$. These illustrative figures confirm that the hybrid block bootstrap has little advantage over MBB in error terms for this problem.

\begin{example}[Nonlinear ARMA(2,3)]
Let $\{X_{t}\}_{t \in \mathbb{Z}}$ be a sequence from the ARMA(2,3) process
\[
X_{t}-0.1X_{(t-1)}+0.3X_{(t-2)}=\epsilon_{t}+0.1\epsilon_{(t-1)}+0.2\epsilon_{(t-2)}-0.1\epsilon_{(t-3)}.
\]
As noted by \citet[Example 6.1]{Lahiri:2003}, such a sequence is strong mixing with exponentially decaying mixing coefficients. To simulate from this model, we initiate by generating $X_{0}, X_{-1}$ from the marginal $N(0, v^{2})$ distribution, which has $v^{2}=1.0776$, with $\epsilon_{0}, \epsilon_{-1}, \epsilon_{-2}$ independent $N(0,1)$. The nonlinear model we consider is the square transformation of the above ARMA process,
\[
Y_{t}=X_{t}^{2}.
\]
\end{example}
The square transformation above preserves the strong mixing property and also preserves the mixing rate. Therefore, $Y_{t}$ is strong mixing with the same exponential rate as $X_{t}$. The interested reader is referred to \citet[p. 69]{FanYao:2003} or \citet[p. 258]{DavisMikosch:2009}. As with the previous example, we consider $p=1/2$, and thus $\xi_{p}$ satisfies
\[
\prob(Y_{t}=X_{t}^{2} \leq \xi_{p})=1/2 ,
\]
implying $\xi_{p}=(0.675v)^{2}$. The simulation approximation to the true value is $G_{n}(-1.5) \approx 0.09276$.

\begin{figure}[h]
\centering
\caption{\label{ARMAsquaredheat}
Heatmap for the nonlinear (squared) ARMA(2,3) model; $n=200$.}
\includegraphics[height=9cm]{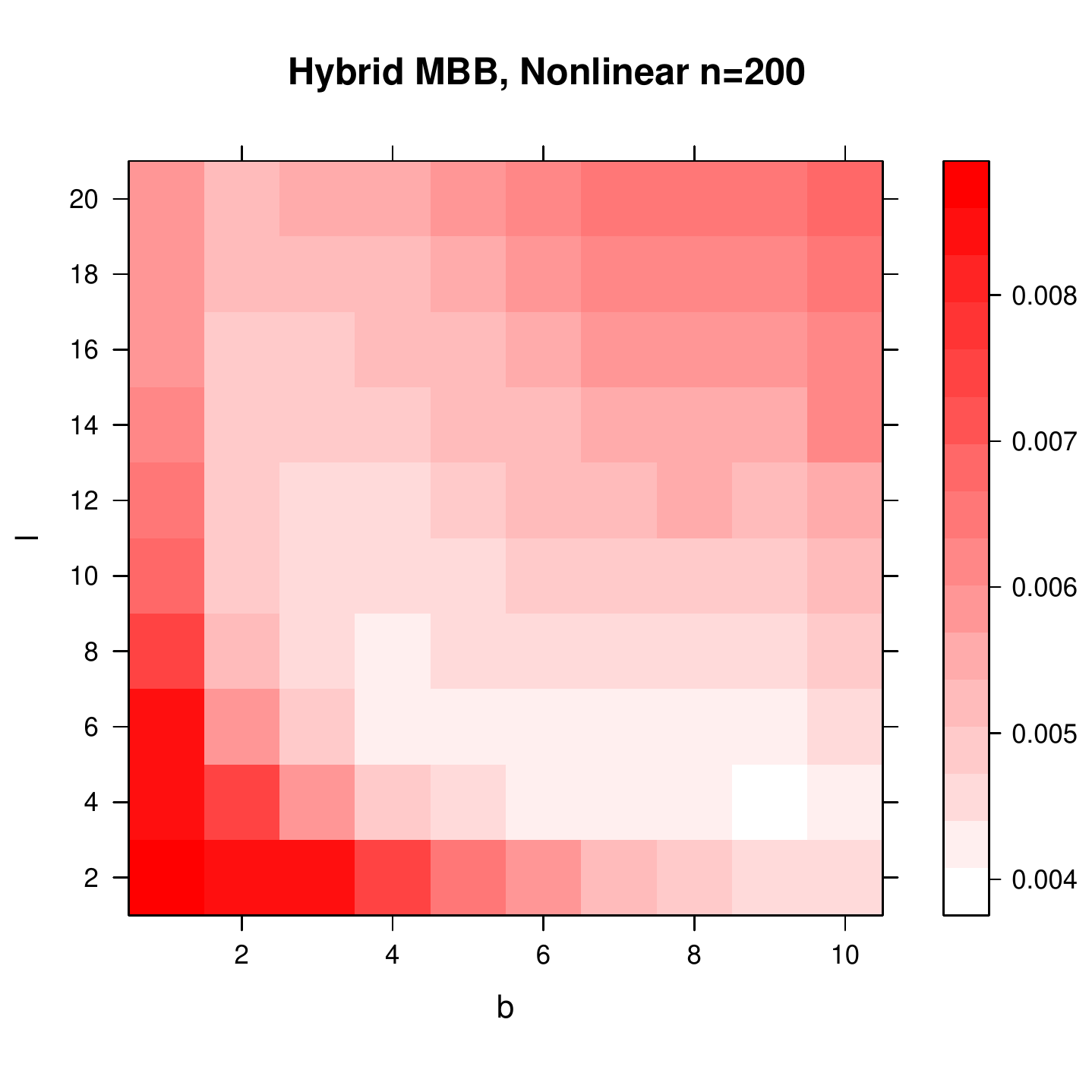}
\end{figure}
The heatmap of Figure~\ref{ARMAsquaredheat} shows again that the subsampling and MBB choices of $(b,\ell)$ are suboptimal from the perspective of minimizing MSE.

\begin{figure}[h]
\centering
\caption{\label{ARMAsquaredheat-lowercoverage}
Heatmap for the coverages of 90\% lower confidence limits in the nonlinear (squared) ARMA(2,3) model; $n=200$.}
\includegraphics[height=9cm]{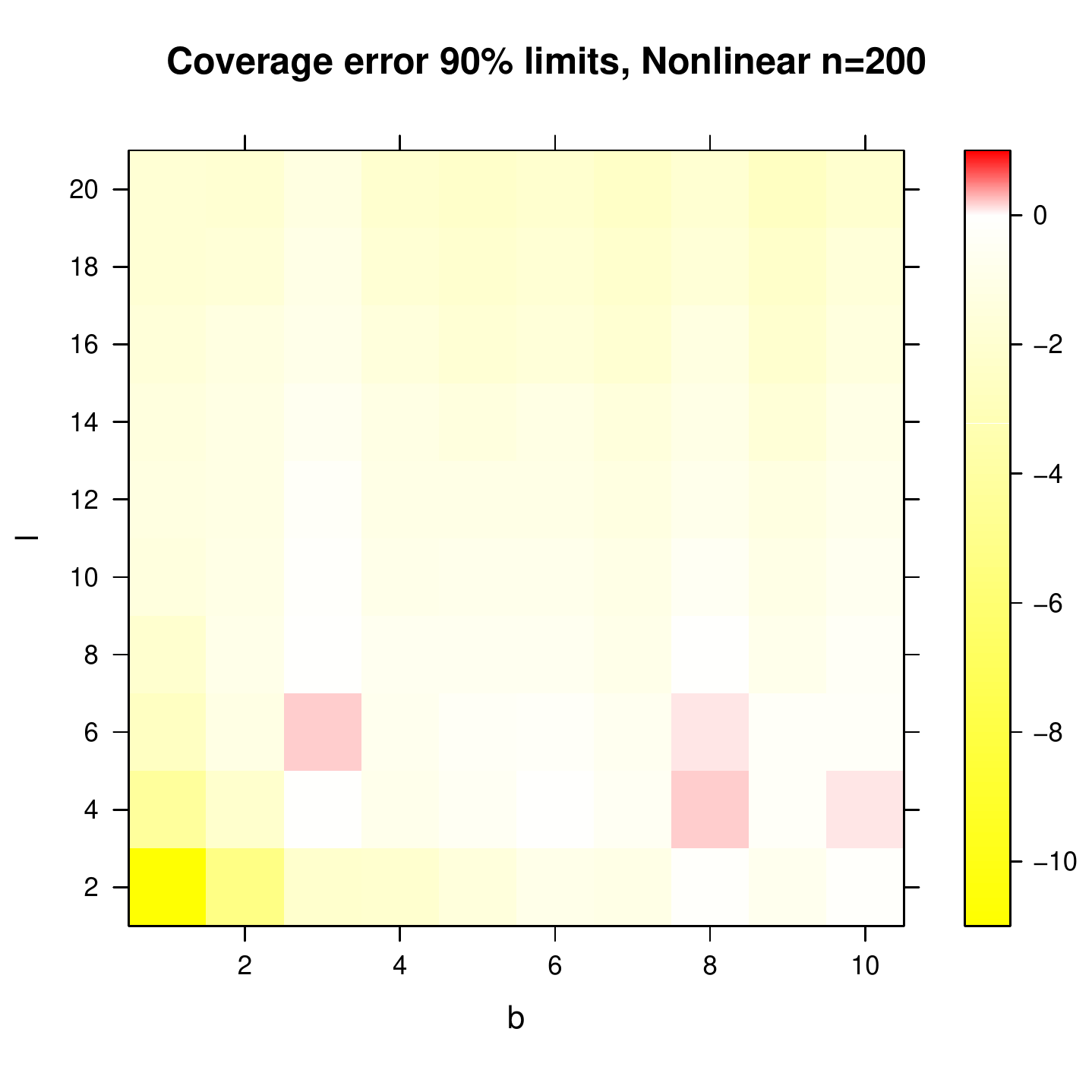}
\end{figure}
In Figure~\ref{ARMAsquaredheat-lowercoverage} we display the coverage error of lower percentile confidence intervals, as described in Section~\ref{cov.error}, of nominal 90\% coverage. We observe that there is undercoverage for most choices of $(b,\ell)$, sometimes very substantial, though there is overcoverage in a few cases. Appropriate choice of $(b, \ell)$ can yield limits with exactly the required coverage.

As proof of concept of the adaptive procedure for choice of $(b, \ell)$ described in Section~\ref{sec:empirical}, we consider estimation of $G_n(1) \approx 0.80952$, for sample size $n=512$. We restrict to candidate values $c_1, c_2 \in \{0.5, 0.75, 1.0, 1.5, 2.0\}$, corresponding to adaptive choice of $b, \ell \in \{4, 6, 8, 12, 16\}$. Table~\ref{tab:adaptive512} shows the MSE in estimation of $G_n(1)$ over 2500 replications for each combination of $(c_1, c_2)$. By contrast, the MSE obtained by minimization of $Err(c_1, c_2)$ for each replication, using 20 subsamples of size $M=64$ in construction of this error quantity, was 0.00189. The adaptive method clearly yields an MSE that is far from optimal in this setting, but outperforms the procedure which fixes $b, \ell$ to larger values among those being considered. 

The adaptive procedure is seen to perform better with increasing sample size. Table~\ref{tab:adaptive1728} provides analagous results for sample size $n=1728$, for which $G_n(1) \approx 0.81125$. Using $M=512$ in the minimization of $Err(c_1, c_2)$ over the same range of $c_1, c_2$, now corresponding to adaptive choice of $b, \ell \in \{6,9,12, 18,24\}$, and again using just 20 subsamples of length $M$ in evaluation of $Err(c_1, c_2)$, the MSE of the adaptively chosen estimator over the 2500 replications was observed as 0.00066, much closer to optimal. Further tuning of the adaptive procedure certainly seems worthwhile as a means of providing an effective automatic choice of $(b, \ell)$ for the hybrid block bootstrap and will be pursued elsewhere.

\begin{table}[h]
\caption{\label{tab:adaptive512}
Nonlinear (squared) ARMA(2,3) model: MSE in estimation of $G_n(1)$ over 2500 replications, for $b= \lfloor c_{1}n^{1/3} \rfloor $ and $\ell= \lfloor c_{2}n^{1/3} \rfloor $, $n=512$.}
\centering
\begin{tabular}{c c | c c c c c}
&&&&$c_2$&&\\ 
&&0.5&0.75&1.0&1.5&2.0 \\
\hline

&0.5&0.00164&0.00150&0.00158&0.00189&0.00230 \\
&0.75&0.00143&0.00154&0.00172&0.00220&0.00272\\
$c_1$&1.0&0.00144&0.00166&0.00191&0.00250& 0.00300\\
&1.5&0.00161&0.00195&0.00232&0.00297& 0.00350\\
&2.0&0.00179&0.00225&0.00265&0.00335& 0.00399\\
\hline

\end{tabular}
\end{table}

\begin{table}[h]
\caption{\label{tab:adaptive1728}
Nonlinear (squared) ARMA(2,3) model: MSE in estimation of $G_n(1)$ over 2500 replications, for $b= \lfloor c_{1}n^{1/3} \rfloor $ and $\ell= \lfloor c_{2}n^{1/3} \rfloor$, $n=1728$.}
\centering
\begin{tabular}{c c | c c c c c}
&&&&$c_2$&&\\ 
&&0.5&0.75&1.0&1.5&2.0 \\
\hline
 &0.5& 0.00062&  0.00063&  0.00072&  0.00089&  0.00108 \\
 &0.75& 0.00061& 0.00131& 0.00082&  0.00105&  0.00126 \\
 $c_1$& 1.0& 0.00065& 0.00080&  0.00094&  0.00119&  0.00139 \\
 &1.5& 0.00076&  0.00094& 0.00112&  0.00138&  0.00166 \\
 &2.0& 0.00087& 0.00106&  0.00125& 0.00159&  0.00182 \\
\hline

\end{tabular}
\end{table}

\begin{figure}[htp]
\centering
\caption{\label{polynomialheat}
Heatmap for the polynomial mixing rate example; $n=200$.}
\includegraphics[height=9cm]{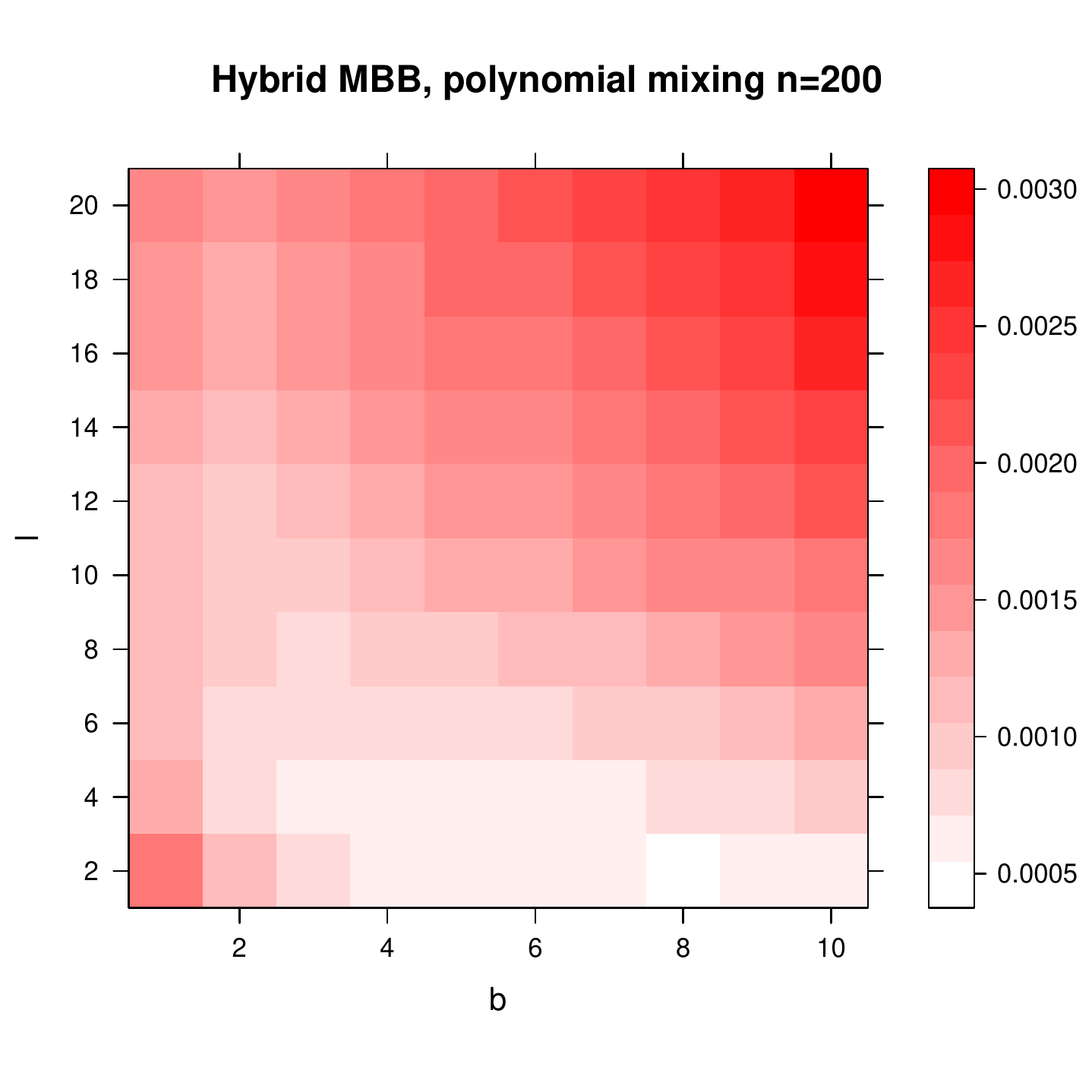}
\end{figure}

Next we construct a process whose mixing coefficients decay at a polynomial rate, but not an exponential rate. This is accomplished through Theorem 2.1 of \citet{Chanda:1974}; see also \citet{Bandyopadhyay:2006} and \S 3 of \citet{CSWX:2016}.
\begin{example}[Polynomial Mixing Rate]
Let the sequence $\{X_{t}\}_{t \in \mathbb{Z}}$ be generated according to
\[
X_{t}=\sum_{j=0}^{\infty}c_{j}Z_{t-j},
\]
where the $Z_{i}$ are independent, identically distributed $N(0,1)$ and $c_{j}=(\tfrac{1}{j+1})^{\nu}$. Then $X_{t}$ is strong mixing with a polynomial rate, and \citet{Chanda:1974} may be used to deduce that $\beta<\nu-2$.

\end{example}
In practice, we cannot simulate from the above process exactly because it is expressed as an infinite series. Therefore, we approximate the process by truncating the series at 100 terms, which means that in reality $X_{t}$ is approximated by a very high order MA process. For our numerical example, $\nu=10.0$ and $n=200$. As before, we consider $p=1/2$, corresponding to $\xi_{p}=0$. Simulation yields an approximation to the true value $G_{n}(2) \approx 0.95229$. The heatmap for this example shown in Figure~\ref{polynomialheat} is based on 10,000 replications, and 10,000 bootstrap samples for each replication of the experiment.
As with the previous two examples, the heatmap supports our finding of suboptimality of the choices of $(b,\ell)$ suggested by the subsampling bootstrap and the MBB.

\section{Proofs}
\label{sec:proofs}
In what follows we denote by $C$ a generic positive constant independent of $n$.
Denote by $\Phi$ and $\phi$ the standard normal distribution and density functions, respectively.
Standard asymptotic properties (e.g. \citet{LahiriSun:2009}) of $F_n$ for dependent data can be invoked to show that, for any $x\in\mathbb{R}$,
\begin{equation}
\prob\Big(\sqrt{n}\big(\hat\xi_n-\xi_p\big)\le x\Big)=
\Phi\big(xf(\xi_p)/\sigma(\xi_p)\big)+O(n^{-1/2}),
\label{pf:true}
\end{equation}
where $\sigma(x)^2=\underset{n\rightarrow\infty}\lim\text{Var}\big(\sqrt{n}F_n(x)\big)
=\sum_{t=-\infty}^\infty\mbox{\rm Cov}\,\big(\pmb{1}\{X_0\le x\},\pmb{1}\{X_t\le x\}\big)$.

We first state a lemma which is a special case of Sun and Lahiri's (2006)
Lemma~5.3.
\begin{lem}
\label{lem:ineq}
Let $\big\{V_{n,t}:t=0,\pm 1,\pm 2,\ldots\big\}$ be a double array of
row-wise stationary strong mixing Bernoulli$\,(p_n)$ random variables with
$0<p_n\le q<1$ and mixing
coefficients $\alpha_n(\cdot)=\alpha(\cdot)$, for some fixed $q\in(0,1)$.
Then, for any positive
$\epsilon_n=o(1)$, $n^{-1}\le\delta_n=o(1)$ and any $\delta\in(0,1)$, we have
\begin{eqnarray*}
\lefteqn{\prob\left(\Big|
\sum_{t=1}^n\big(V_{n,t}-p_n\big)\Big|>n\epsilon_n\right)}\nonumber\\
&\le& C
\left(\delta_n^{-1}+\dfrac{\epsilon_n^2}{p_n+\epsilon_n}\right)
\exp\left\{-\,\dfrac{Cn\delta_n\epsilon_n^2}{p_n+\epsilon_n}\right\}
+Cn\big(1+p_n^{\delta}\epsilon_n^{-1}\big)\delta_n^{\beta(1-\delta)}.
\end{eqnarray*}
\end{lem}
Define, for any $r>0$, $\mathscr{B}_r(\xi_p)=[\xi_p-r,\xi_p+r]$.
\begin{lem}
\label{lem:empcdf}
Suppose that $n^{-\frac{4\beta+7}{6(3\beta+5)}}\ell\rightarrow\infty$. Then
for any arbitrarily small $\delta>0$,
the following results hold uniformly over $\epsilon\in\big[n^{-c_0},1\big)$.
\begin{itemize}
\item[(i)]
$\displaystyle\underset{x\in\mathscr{B}_\epsilon(\xi_p)\cap\mathscr{N}_p}\sup
\big|F_n(x)-F(x)\big|=O_p\left(n^{-\,\frac{\beta-1}{2(\beta+1)}+3\delta}
\epsilon^{\frac{1}{2(\beta+1)}+\delta}\right)$
for any $c_0\in(0,3)$.

\item[(ii)]
$\displaystyle\underset{x\in\mathscr{B}_\epsilon(\xi_p)\cap\mathscr{N}_p}\sup
\big|\tilde{F}_n(x)-F_n(x)\big|=O_p\left(
n^{-1}\epsilon^{\frac{1}{2(\beta+1)}+\delta}
\ell^{\frac{\beta+3}{2(\beta+1)}+3\delta}\right)$
for some $c_0>1/2$.

\item[(iii)]
$\displaystyle\underset{x\in\mathscr{B}_\epsilon(\xi_p)\cap\mathscr{N}_p}\sup
\big|{F}_n(x)-F_n(\xi_p)-F(x)+p\big|=O_p\left(n^{-\frac{\beta-1}{2(\beta+1)}+\delta}\epsilon^{(1+\delta)/2}\right)$
for any $c_0\in(0,2)$.
\end{itemize}
\end{lem}
{\it Proof of Lemma~\ref{lem:empcdf}\/}:

Define $m_n=\left\lceil n^{\,\frac{\beta-1}{2(\beta+1)}-3\delta}
\epsilon^{1-\frac{1}{2(\beta+1)}-\delta}\right\rceil$ and $\epsilon_n=\epsilon/m_n$.
Then we have
\begin{equation}
\label{pf:empcdf1}
\underset{x\in\mathscr{B}_\epsilon(\xi_p)}\sup
\big|F_n(x)-F(x)\big|
\le \underset{k\in\{0,\pm 1,\ldots,\pm m_n\}}
\max\big|F_n(\xi_p+k\epsilon_n)-F(\xi_p+k\epsilon_n)\big|+C\epsilon_n.
\end{equation}
For each $k\in\{0,\pm 1,\ldots,\pm m_n\}$, application of
Lemma~\ref{lem:ineq} with
$V_{n,t}=\pmb{1}\{X_t\le\xi_p+k\epsilon_n\}$,
$p_n=F(\xi_p+k\epsilon_n)$ and
$\delta_n=\big(\epsilon n^{2-\Delta_1}\big)^{\frac{-1}{1+\beta(1-\Delta_2)}}$,
for arbitrarily small $\Delta_1,\Delta_2>0$, yields
\begin{eqnarray*}
\lefteqn{\prob\left(\Big|
F_n(\xi_p+k\epsilon_n)-F(\xi_p+k\epsilon_n)\Big|>\epsilon_n\right)}\nonumber\\
&\le&C
n^{\frac{4}{2+\beta}}
\exp\left\{-C
n^{(3-c_0)\delta}
\right\}+C\epsilon^{-\frac{1 + 3 \beta + 2 \beta^2}{2 (1 + \beta)^2}-
\frac{\delta}{
\beta+1} - C\Delta_2}
n^{-\frac{\beta - 1}{2 ( \beta + 1)} -
\frac{(\beta+ 3)\delta}{\beta+1}+C\Delta_2}.
\end{eqnarray*}
It follows by Bonferroni's inequality that
\begin{eqnarray*}
\lefteqn{\prob\left(\underset{k\in\{0,\pm 1,\ldots,\pm m_n\}}
\max\big|F_n(\xi_p+k\epsilon_n)-F(\xi_p+k\epsilon_n)\big|>\epsilon_n\right)}\nonumber\\
&\le&(2m_n+1)
\prob\left(\Big|
F_n(\xi_p+k\epsilon_n)-F(\xi_p+k\epsilon_n)\Big|>\epsilon_n\right)\\
&\le&Cn^{\frac{4}{2+\beta}+\frac{\beta-1}{2(\beta+1)}-3\delta}
\exp\left\{-C
n^{(3-c_0)\delta}
\right\}+C
n^{-\left(
\frac{(4-c_0)\beta+ 6-2c_0}{\beta+1}\right)\delta+C\Delta_2}=o(1)
\end{eqnarray*}
for sufficiently small $\Delta_2$,
uniformly over $\epsilon\in\big[n^{-c_0},1\big)$.
This, in conjunction with (\ref{pf:empcdf1}), implies  that
$\underset{x\in\mathscr{B}_\epsilon(\xi_p)}\sup
\big|F_n(x)-F(x)\big|=O_p(\epsilon_n)$, which proves part (i) of
Lemma~\ref{lem:empcdf}.

To prove part (ii), write $n'=n-\ell+1$ and note that
\begin{eqnarray}
F_n(x)-\tilde{F}_n(x)&=&
\dfrac{1}{\ell n'}\sum_{i=1}^{\ell-1}
\left(\dfrac{\ell n'}{n}-i\right)\big(\pmb{1}\{X_i\le x\}
+\pmb{1}\{X_{n+1-i}\le x\}-2F(x)\big)\nonumber\\
&&\qquad-\,\dfrac{\ell-1}{nn'}\sum_{i=\ell}^{n'}\big(\pmb{1}\{X_{i}\le x\}-F(x)\big).
\label{pf:empcdf2}
\end{eqnarray}
Define, for $c_0>1/2$, $\displaystyle m=
n^{-\frac{c_0}{4(1+\beta)}-\frac{c_0\delta}{2}}\ell^{\frac{3\beta+5}{4(\beta+1)}+\frac{3\delta}{2}}$.
It is clear that $m=o(\ell)$. Noting that for sufficiently large $n$
and sufficiently small $\Delta>0$,
\[ \dfrac{\ln\ell}{\ln n}
>\frac{4\beta+7}{6(3\beta+5)}+\Delta, \]
we have
\begin{eqnarray*}
\dfrac{\ln m}{\ln n}
&>&-c_0\left(\dfrac{1}{4(\beta+1)}+\dfrac{\delta}{2}\right)+
\left(\dfrac{3\beta+5}{4(\beta+1)}+\dfrac{3\delta}{2}\right)\left(\frac{4\beta+7}{6(3\beta+5)}+\Delta\right)\\
&=&-(c_0-1/2)\left(\dfrac{1}{4(\beta+1)}+\dfrac{\delta}{2}\right)+\dfrac{1}{6}+
\dfrac{(3\beta+5)\Delta}{4(\beta+1)}+\dfrac{\delta}{4}\left(\frac{\beta+2}{3\beta+5}+6\Delta\right).
\end{eqnarray*}
We may therefore choose $c_0$ sufficiently close to $1/2$ and some $K<6$ such that
$\ln n/\ln m\le K$ and $c_0K<3$. It follows that for any $\epsilon\ge n^{-c_0}$,
\[ \epsilon\ge n^{-c_0}\ge m^{-c_0K}>\ell^{-c_0K}. \]
Consider
\begin{eqnarray}
\lefteqn{
\dfrac{1}{\ell n'}\sum_{i=1}^{\ell-1}
\left(\dfrac{\ell n'}{n}-i\right)\big(\pmb{1}\{X_{i}\le x\}-F(x)\big)}\nonumber\\
&=&
\dfrac{1}{n}\sum_{i=1}^{\ell-1}\big(\pmb{1}\{X_i\le x\}-F(x)\big)-\dfrac{1}{\ell n'}
\sum_{j=1}^{\ell-m}\left\{\sum_{i=j}^{\ell-1}\big(\pmb{1}\{X_i\le x\}-F(x)\big)\right\}\nonumber\\
&&-\,\dfrac{1}{\ell n'}
\sum_{j=\ell-m+1}^{\ell-1}\left\{\sum_{i=j}^{\ell-1}\big(\pmb{1}\{X_i\le x\}-F(x)\big)\right\}
=I_1-I_2-I_3, \text{\ say.}\nonumber
\end{eqnarray}
Applying part (i), we have, uniformly over $\epsilon\in\big[n^{-c_0},1\big)$, that
\[
I_1=O_p\left(n^{-1}\ell^{\frac{\beta+3}{2(\beta+1)}+3\delta}
\epsilon^{\frac{1}{2(\beta+1)}+\delta}\right)
\]
and
\[
I_2=O_p\left(
\dfrac{1}{\ell n'}
\sum_{j=1}^{\ell-m}(\ell-j)^{\frac{\beta+3}{2(\beta+1)}+3\delta}
\epsilon^{\frac{1}{2(\beta+1)}+\delta}\right)=O_p\left(n^{-1}\ell^{\frac{\beta+3}{2(\beta+1)}+3\delta}
\epsilon^{\frac{1}{2(\beta+1)}+\delta}\right).
\]
It is clear that
\begin{align*}
I_3&=O_p\left((n\ell)^{-1}m^2\right) \\
&=O_p\left(
n^{-1-\frac{c_0}{2(1+\beta)}-c_0\delta}\ell^{\frac{\beta+3}{2(\beta+1)}+3\delta}\right) \\
&=O_p\left(n^{-1}\ell^{\frac{\beta+3}{2(\beta+1)}+3\delta}
\epsilon^{\frac{1}{2(\beta+1)}+\delta}\right).
\end{align*}
The bounds on $I_1,I_2,I_3$ therefore imply that
\begin{equation}
\begin{split}
\dfrac{1}{\ell n'}\sum_{i=1}^{\ell-1}
\left(\dfrac{\ell n'}{n}-i\right)\big(\pmb{1}\{X_i\le x\}
+\pmb{1}\{X_{n+1-i}\le x\}-2F(x)\big)\\
\qquad=O_p\left(n^{-1}\ell^{\frac{\beta+3}{2(\beta+1)}+3\delta}
\epsilon^{\frac{1}{2(\beta+1)}+\delta}\right).
\end{split}
\label{pf:empcdf3}
\end{equation}
It follows by part (i) again that
\begin{eqnarray}
\dfrac{\ell-1}{nn'}\sum_{i=\ell}^{n'}\big(\pmb{1}\{X_{i}\le x\}-F(x)\big)
&=&O_p\left(\ell n^{-1-\,\frac{\beta-1}{2(\beta+1)}+3\delta}
\epsilon^{\frac{1}{2(\beta+1)}+\delta}\right)\nonumber\\
&=&o_p\left(n^{-1}\ell^{\frac{\beta+3}{2(\beta+1)}+3\delta}
\epsilon^{\frac{1}{2(\beta+1)}+\delta}\right).
\label{pf:empcdf4}
\end{eqnarray}
Part (ii) then follows by combining (\ref{pf:empcdf2}), (\ref{pf:empcdf3}) and (\ref{pf:empcdf4}).

For the proof of part (iii), define $M_n=\left\lceil n^{\frac{\beta-1}{2(\beta+1)}-\delta}\epsilon^{(1-\delta)/2}\right\rceil$ and $\epsilon_n=\epsilon/M_n$.
Then we have
\begin{eqnarray}
\lefteqn{
\underset{x\in\mathscr{B}_\epsilon(\xi_p)}\sup
\big|F_n(x)-F_n(\xi_p)-F(x)+p\big|}\nonumber\\
&\le& \underset{k\in\{0,\pm 1,\ldots,\pm M_n\}}
\max\big|F_n(\xi_p+k\epsilon_n)-F_n(\xi_p)-F(\xi_p+k\epsilon_n)+p\big|+C\epsilon_n.
\label{pf:empcdf5}
\end{eqnarray}
For each $k\in\{\pm 1,\ldots,\pm m_n\}$, set, for the application of
Lemma~\ref{lem:ineq},
$V_{n,t}=\big|\pmb{1}\{X_t\le\xi_p+k\epsilon_n\}-\pmb{1}\{X_t\le\xi_p\}\big|$,
$p_n=\big|F(\xi_p+k\epsilon_n)-p\big|$ and
$\delta_n=\big(n^{2-\Delta_1}\epsilon^{\Delta_2}\big)^{-1/(\beta(1-\Delta_2)+1)}$,
for arbitrarily small $\Delta_1,\Delta_2>0$. Noting that
$C^{-1}\epsilon_n\le p_n\le C\epsilon$, we have, by Lemma~\ref{lem:ineq},
that for any $\Delta\in(0,1)$,
\begin{eqnarray*}
\lefteqn{\prob\left(\Big|
F_n(\xi_p+k\epsilon_n)-F_n(\xi_p)-F(\xi_p+k\epsilon_n)+p\Big|>\epsilon_n\right)}\nonumber\\
&\le& C
\big(n^{2-\Delta_1}\epsilon^{\Delta_2}\big)^{1/(\beta(1-\Delta_2)+1)}
\exp\left\{-\,C\big(n^{2-\Delta_1}\epsilon^{\Delta_2}\big)^{-1/(\beta(1-\Delta_2)+1)}n^{\frac{2}{\beta+1}+2\delta}\epsilon^{\delta}\right\}\\
&&+\,Cn^{\frac{3\beta+1}{2(\beta+1)}-\delta}\epsilon^{-(1+\delta)/2+\Delta}
\big(n^{2-\Delta_1}\epsilon^{\Delta_2}\big)^{-\beta(1-\Delta)/(\beta(1-\Delta_2)+1)}.
\end{eqnarray*}
It follows by Bonferroni's inequality that for sufficiently small $\Delta'>0$ and
for any $c_0\in(0,2)$,
\begin{eqnarray*}
\lefteqn{\prob\left(\underset{k\in\{0,\pm 1,\ldots,\pm m_n\}}
\max\big|F_n(\xi_p+k\epsilon_n)-F_n(\xi_0)-
F(\xi_p+k\epsilon_n)+p\big|>\epsilon_n\right)}\nonumber\\
&\le&(2M_n+1)
\prob\left(\Big|
F_n(\xi_p+k\epsilon_n)-F_n(\xi_0)-
F(\xi_p+k\epsilon_n)+p\Big|>\epsilon_n\right)\\
&\le&
Cn^{\frac{\beta+3}{2 (\beta+1)}}
\exp\left\{-\,Cn^{(2-c_0)\delta/2}
\right\}+Cn^{-(2-c_0)\delta+\Delta'}=o(1),
\end{eqnarray*}
uniformly over $\epsilon\in\big[n^{-c_0},1\big)$.
This, in conjunction with (\ref{pf:empcdf5}), implies  that
$\underset{x\in\mathscr{B}_\epsilon(\xi_p)}\sup
\big|F_n(x)-F_n(\xi_p)-F(x)+p\big|=O_p(\epsilon_n)$, which proves part (iii) of
Lemma~\ref{lem:empcdf}.

\begin{lem}
\label{lem:quan}
Suppose that $n^{-\frac{4\beta+7}{6(3\beta+5)}}\ell\rightarrow\infty$. Then
for any arbitrarily small $\delta>0$,
\begin{itemize}
\item[(i)] $\displaystyle
\tilde\xi_n=\xi_p+O_p\left(n^{-1/2}+n^{-\frac{2(\beta+1)}{2\beta+1}+\delta}
\ell^{\frac{\beta+3}{2\beta+1}}\right)$.
\item[(ii)] $\displaystyle
\tilde{F}_n\big(\tilde\xi_n\big)=p+o_p\Big(
n^{-\frac{\beta-3}{\beta-1}+\delta}
+n^{-\frac{\beta(2\beta-3)}{(\beta-1)(2\beta+1)}+\delta}
\ell^{\frac{2(\beta+3)}{(\beta-1)(2\beta+1)}}
+n^{-\frac{4\beta+5}{4(\beta+1)}+\delta}
\ell^{\frac{\beta+3}{2(\beta+1)}}$

\qquad\qquad\qquad\qquad$+\,n^{-\frac{2(\beta+1)}{2\beta+1}+\delta}
\ell^{\frac{\beta+3}{2\beta+1}}
\Big)$.
\end{itemize}
\end{lem}
{\it Proof of Lemma~\ref{lem:quan}\/}:

Let $c_0>1/2$ be as specified in Lemma~\ref{lem:empcdf}(ii). Define, for $\epsilon\in[n^{-c_0},1)$,
\[\delta_n(\epsilon)=
n^{-1}\epsilon^{\frac{1}{2(\beta+1)}+\delta}
\ell^{\frac{\beta+3}{2(\beta+1)}+3\delta},
\]
\[
\epsilon_1=\left\{n^{-1}
\ell^{\frac{\beta+3}{2(\beta+1)}+3\delta}\right\}^{
\frac{2(\beta+1)}{2\beta+1-2\delta(\beta+1)}}
\mbox{\ \ and\ \ }
\epsilon_2=n^{-1/2}+\epsilon_1.
\]
Note that $\delta_n(\epsilon_1)=\epsilon_1$.
Using Lemma~\ref{lem:empcdf}(ii), we have, for some $\alpha\in(0,1)$, any $M>0$ and sufficiently large $n,\tilde{M}$,
\begin{eqnarray}
\lefteqn{\prob\big(\tilde\xi_n-\xi_p>M\epsilon_2\big)}\nonumber\\
&\le& \prob\Big(\underset{x\in\mathscr{B}_{M\epsilon_2}(\xi_p)\cap\mathscr{N}_p}\sup
|\tilde{F}_n(x)-F_n(x)|>\tilde{M}\delta_n(M\epsilon_2)\Big)\nonumber\\
&&\qquad+\,\prob\big(F_n(\xi_p+M\epsilon_2)\le p+\tilde{M}\delta_n(M\epsilon_2)\big)\nonumber\\
&\le&\tilde{M}^{-1}+\prob\Big(\sqrt{n}\big(F_n(\xi_p+M\epsilon_2)-F(\xi_p+M\epsilon_2)\big)\nonumber\\
&&\;\;\le -2^{-1}Mf(\xi_p)\big\{1+
\sqrt{n}\epsilon_1-CM^{-(1-\alpha)}\tilde{M}\sqrt{n}\big(\delta_n(n^{-1/2})+\epsilon_1\big)\big\}\Big).
\label{pf:quan1}
\end{eqnarray}
Note that if $\epsilon_1=o(n^{-1/2})$, then $\delta_n(n^{-1/2})=o(n^{-1/2})$, so that
\[ \underset{n\rightarrow\infty}\lim
-2^{-1}Mf(\xi_p)\big\{1+
\sqrt{n}\epsilon_1-CM^{-(1-\alpha)}\tilde{M}\sqrt{n}\big(\delta_n(n^{-1/2})+\epsilon_1\big)\big\}
=-2^{-1}Mf(\xi_p).\]
If $\epsilon_1\ge C_0n^{-1/2}$ for some $C_0>0$, then $\delta_n(n^{-1/2})\le C_0^{-\alpha}\delta_n(\epsilon_1)=C_0^{-\alpha}\epsilon_1$,
so that for $M$ sufficiently large, we have
\begin{eqnarray*}
\lefteqn{\sqrt{n}\epsilon_1-CM^{-(1-\alpha)}\tilde{M}\sqrt{n}\big(\delta_n(n^{-1/2})+\epsilon_1\big)}\\
&\ge& C_0\left\{1-CM^{-(1-\alpha)}\tilde{M}\big(C_0^{-\alpha}+1\big)\right\}\ge C_0/2.
\end{eqnarray*}
The above results suggest that
\[ -2^{-1}Mf(\xi_p)\big\{1+
\sqrt{n}\epsilon_1-CM^{-(1-\alpha)}\tilde{M}\sqrt{n}\big(\delta_n(n^{-1/2})+\epsilon_1\big)\big\}\le
-4^{-1}Mf(\xi_p) \]
for $M,n$ sufficiently large. It then follows from (\ref{pf:quan1}) that
$\prob\big(\tilde\xi_n-\xi_p>M\epsilon_2\big)$ can be made arbitrarily small if we choose $n,M$ and $\tilde{M}$
large enough, using the fact that $F_n(\xi_p+M\epsilon_2)-F(\xi_p+M\epsilon_2)=O_p(n^{-1/2})$. The same arguments can be applied to
the lower tail probability $\prob\big(\tilde\xi_n-\xi_p<-M\epsilon_2\big)$.
Thus we have $\tilde\xi_n=\xi_p+O_p(\epsilon_2)$. Lemma~\ref{lem:quan}(i) then
follows as $\delta$ can be made arbitrarily small.

To prove part (ii), write
\begin{align*}
\varepsilon_0(n,\ell)&=n^{-1/2}+n^{-\frac{2(\beta+1)}{2\beta+1}}
\ell^{\frac{\beta+3}{2\beta+1}}, \\
\varepsilon_1(n,\ell)&=n^{-\frac{4\beta+5}{4(\beta+1)}}
\ell^{\frac{\beta+3}{2(\beta+1)}}+n^{-\frac{2(\beta+1)}{2\beta+1}}
\ell^{\frac{\beta+3}{2\beta+1}},\\
\varepsilon_2(n,\ell)&=n^{-\frac{\beta-3}{\beta-1}}
+n^{-\frac{\beta(2\beta-3)}{(\beta-1)(2\beta+1)}}
\ell^{\frac{2(\beta+3)}{(\beta-1)(2\beta+1)}}, \\
\varepsilon(n,\ell)&=\varepsilon_1(n,\ell)+\varepsilon_2(n,\ell).
\end{align*}
Denote by $X_{(1)}\le\cdots\le X_{(n)}$ the ordered sequence of $X_1,\ldots,X_n$.
For any arbitrarily small $\Delta,\Delta'>0$, $\prob\big(\tilde{F}_n(\tilde\xi_n)\ge p+n^\delta\varepsilon(n,\ell)\big)$ is bounded above by
\begin{equation}
\begin{split}
\prob\Big(F_n(X_{(j+1)})=\cdots=F_n(X_{(j+k)})\ge p+n^\delta\varepsilon(n,\ell)
-n^\Delta\varepsilon_1(n,\ell),\qquad\qquad\qquad\\
F_n(X_{(j)})<p+n^\Delta\varepsilon_1(n,\ell),\;
\big|X_{(j+1)}-\xi_p\big|\le n^{\Delta'}\varepsilon_0(n,\ell)\text{\ some\ }
j,k\ge1 \Big)\\
+\,\prob\Big(\big|\tilde\xi_n-\xi_p\big|> n^{\Delta'}\varepsilon_0(n,\ell)\Big)+\prob\Big(\sup
\big|\tilde{F}_n(x)-F_n(x)\big|>n^\Delta\varepsilon_1(n,\ell)\Big),
\end{split}
\label{pf:quan2}
\end{equation}
where the supremum in the last probability is taken over
$x\in\mathscr{B}_{n^{\Delta'}\varepsilon_0(n,\ell)}(\xi_p)
\cap\mathscr{N}_p$.
Noting that $n^{4-\beta}\varepsilon_0(n,\ell)^2=O\big(\varepsilon_2(n,\ell)^{\beta-1}\big)=O\big(\varepsilon(n,\ell)^{\beta-1}\big)$,
we have
\begin{eqnarray*}
\lefteqn{n^{-(\beta-2-\Delta')}\varepsilon_0(n,\ell)
\big(n^\delta\varepsilon(n,\ell)\big)^{-
(\beta-1)}} \\
&=&O\Big(n^{-(\beta-2-\Delta')-\delta(\beta-1)}\varepsilon_0(n,\ell)
n^{\beta-4}\varepsilon_0(n,\ell)^{-2}\Big)\\
&=&O\Big(n^{-2+\Delta'-\delta(\beta-1)}\varepsilon_0(n,\ell)^{-1}
\Big)\\
&=&o\Big(
n^{-2-\Delta'}\varepsilon_0(n,\ell)^{-1}\Big)
\end{eqnarray*}
for sufficiently small $\Delta'>0$. Thus we may find a positive sequence $\{\eta_n\}$ satisfying
\begin{equation}
\eta_n=o\big(n^{-2-\Delta'}\varepsilon_0(n,\ell)^{-1}\big)\mbox{\ \ and\ \ }
n^{-(\beta-2-\Delta')}\varepsilon_0(n,\ell)
\big(n^\delta\varepsilon(n,\ell)\big)^{-
(\beta-1)}=o(\eta_n)
\label{pf:quan3}
\end{equation}
for sufficiently small $\Delta'>0$.
Noting that $n^\Delta\varepsilon_1(n,\ell)=o\big(n^\delta\varepsilon(n,\ell)\big)$
for sufficiently small $\Delta$, following the proof of Lemma~5.4(iv) of Sun and Lahiri (2006)
and using (\ref{pf:quan3}),
the first probability  in (\ref{pf:quan2})
can be bounded above, for sufficiently large $n$ and sufficiently small $\Delta'>0$, by
\begin{eqnarray}
\lefteqn{n\hspace{-3ex}\sum_{j>2^{-1}n^{1+\delta}\varepsilon(n,\ell)}
\hspace{-2ex}\prob\Big(X_1=X_j\in\mathscr{B}_{n^{\Delta'}\varepsilon_0(n,\ell)}(\xi_p)
\cap\mathscr{N}_p\Big)}\nonumber\\
&\le&Cn^{2+\Delta'}\varepsilon_0(n,\ell)\eta_n+Cn^{2-\beta+\Delta'}\varepsilon_0(n,\ell)
\big(n^{\delta}\varepsilon(n,\ell)\big)^{-(\beta-1)}\eta_n^{-1}=o(1).
\nonumber
\end{eqnarray}
That the last two probabilities in (\ref{pf:quan2})
converge to 0 follows
from Lemma~\ref{lem:quan}(i) and Lemma~\ref{lem:empcdf}(ii), respectively.
From the above results we derive that $\prob\big(\tilde{F}_n(\tilde\xi_n)\ge p+n^\delta\varepsilon(n,\ell)\big)=o(1)$.
Similar arguments show also that $\prob\big(\tilde{F}_n(\tilde\xi_n)\le p-n^\delta\varepsilon(n,\ell)\big)=o(1)$, which
completes the proof of part (ii).

\begin{lem}
\label{lem:bootcdf}
For any arbitrarily small $\delta>0$ and any compact $\mathscr{K}\subset\mathbb{R}$,
\begin{eqnarray*}
\lefteqn{\prob\Big(\sqrt{b\ell}\big(F^*_n(x)-\tilde{F}_n(x)\big)\le y\Big|X_1,\ldots,X_n\Big)}\\
&=&\Phi\big(y/\sigma(x)\big)
+
\begin{cases}
O_p\big(\ell^{-1}+\ell^{1/2}n^{-1/2}+(b\ell)^{-1/2}\ell^\delta\big),&
\ell=O(b),\\
O_p\big(\ell^{-1}+\ell^{1/2}n^{-1/2}+(b\ell)^{-1/2}\big),&
\beta=\infty,
\end{cases}
\end{eqnarray*}
uniformly over $(x,y)\in\mathscr{N}_p\times\mathscr{K}$.
\end{lem}
{\it Proof of Lemma~\ref{lem:bootcdf}\/}:

Denote by $\hat\kappa_j(x)$ the $j$th conditional cumulant of
$\sqrt{b\ell}\big(F^*_n(x)-\tilde{F}_n(x)\big)$ given $X_1,\ldots,X_n$. It is clear that
$\hat\kappa_1(x)=0$. It can be shown using somewhat tedious asymptotic
 expansions that, uniformly over $x\in\mathscr{N}_p$,
\begin{align*}
\hat\kappa_2(x)&=\sigma(x)^2+O_p\big(\ell^{-1}+\ell^{1/2}n^{-1/2}\big),\\
\hat\kappa_j(x)&=O_p\Big(b^{-(j-2)/2}
\big\{\ell^{-1/2+\delta(\beta)}+\ell^{1/2}n^{-1/2}+\ell^{(j+1-\beta)/2}n^{-1/2}\big\}\Big),
\;\;\;j\ge 3,
\end{align*}
where $\delta(\beta)>0$ can be made arbitrarily small if $\beta<\infty$ and $\delta(\infty)=0$.
Lemma~\ref{lem:bootcdf} then follows by comparing the Edgeworth expansion of
the conditional distribution function of $\sqrt{b\ell}\big(F^*_n(x)-\tilde{F}_n(x)\big)$ with
the $N\big(0,\sigma(x)^2\big)$ distribution function.

\noindent
{\it Proof of Theorem~\ref{thm:main}\/}:

Note first that, by Lemmas~\ref{lem:empcdf}, \ref{lem:quan} and Taylor expansion of $F$ about
$\tilde\xi_n$,
\begin{eqnarray}
\lefteqn{p-\tilde{F}_n\big(\tilde\xi_n+(b\ell)^{-1/2}x\big)}\nonumber\\
&=&\Big\{p-\tilde{F}_n(\tilde\xi_n)\Big\}+
\Big\{\tilde{F}_n(\tilde\xi_n)-\tilde{F}_n\big(\tilde\xi_n+(b\ell)^{-1/2}x\big)\Big\}\nonumber\\
&=& {F}_n(\tilde\xi_n)-{F}_n\big(\tilde\xi_n+(b\ell)^{-1/2}x\big)+
o_p\Big(
n^{-\frac{\beta-3}{\beta-1}+\delta}
+n^{-\frac{\beta(2\beta-3)}{(\beta-1)(2\beta+1)}+\delta}
\ell^{\frac{2(\beta+3)}{(\beta-1)(2\beta+1)}}\nonumber\\
&&+\,n^{-\frac{4\beta+5}{4(\beta+1)}+\delta}
\ell^{\frac{\beta+3}{2(\beta+1)}}+n^{-\frac{2(\beta+1)}{2\beta+1}+\delta}
\ell^{\frac{\beta+3}{2\beta+1}}
\Big)+O_p\Big(n^{-1}b^{-\frac{1}{4(\beta+1)}-\delta}
\ell^{\frac{2\beta+5}{4(\beta+1)}+5\delta}\Big)\nonumber\\
&=&-(b\ell)^{-1/2}xf(\tilde\xi_n)+
o_p\Big(
n^{-\frac{\beta-3}{\beta-1}+\delta}+n^{-\frac{3\beta-1}{4(\beta+1)}+\delta}
+n^{-\frac{\beta(2\beta-3)}{(\beta-1)(2\beta+1)}+\delta}
\ell^{\frac{2(\beta+3)}{(\beta-1)(2\beta+1)}}\nonumber\\
&&+\,n^{-\frac{4\beta+5}{4(\beta+1)}+\delta}
\ell^{\frac{\beta+3}{2(\beta+1)}}+n^{-\frac{2(\beta+1)}{2\beta+1}+\delta}
\ell^{\frac{\beta+3}{2\beta+1}}
+n^{-\frac{4\beta^2+3\beta+1}{2(2\beta+1)(\beta+1)}+\delta}\ell^{\frac{\beta+3}{2(2\beta+1)}}\Big)
\nonumber\\
&& \;\; +\,O_p\left((b\ell)^{-1}+
n^{-1}b^{-\frac{1}{4(\beta+1)}-\delta}
\ell^{\frac{2\beta+5}{4(\beta+1)}+5\delta}
+n^{-\frac{\beta-1}{2(\beta+1)}+\delta}(b\ell)^{-(1+\delta)/4}\right).
\label{pf:main1}
\end{eqnarray}
Applying Lemma~\ref{lem:bootcdf}, we have, for arbitrarily small $\delta>0$, that
\begin{eqnarray}
\lefteqn{\prob\Big(F^*_n\big(\tilde\xi_n+(b\ell)^{-1/2}x\big)\le p\Big|X_1,\ldots,X_n\Big)}\nonumber\\
&=&\prob\Big(\sqrt{b\ell}\big\{F^*_n\big(\tilde\xi_n+(b\ell)^{-1/2}x\big)
-\tilde{F}_n\big(\tilde\xi_n+(b\ell)^{-1/2}x\big)\big\}\nonumber\\
&&\qquad\le \sqrt{b\ell}\big\{p-\tilde{F}_n\big(\tilde\xi_n+(b\ell)^{-1/2}x\big)\big\}
\Big|X_1,\ldots,X_n\Big)\nonumber\\
&=&\Phi\Big(\sqrt{b\ell}\big\{p-\tilde{F}_n\big(\tilde\xi_n+(b\ell)^{-1/2}x\big)\big\}
/\sigma\big(\tilde\xi_n+(b\ell)^{-1/2}x\big)\Big)\nonumber\\
&&+\begin{cases}
O_p\big(\ell^{-1}+\ell^{1/2}n^{-1/2}+(b\ell)^{-1/2}\ell^\delta\big),&
\ell=O(b),\\
O_p\big(\ell^{-1}+\ell^{1/2}n^{-1/2}+(b\ell)^{-1/2}\big),&
\beta=\infty.
\end{cases}
\label{pf:main2}
\end{eqnarray}
It follows from (\ref{pf:main1}), (\ref{pf:main2}) and Lemma~\ref{lem:quan}(i) that for arbitrarily small $\delta>0$,
\begin{eqnarray}
\lefteqn{\prob\Big(F^*_n\big(\tilde\xi_n+(b\ell)^{-1/2}x\big)\le p\Big|X_1,\ldots,X_n\Big)}\nonumber\\
&=&\Phi\big(-xf(\xi_p)
/\sigma(\xi_p)\big)+
O_p\Big(\ell^{-1}+\ell^{1/2}n^{-1/2}+(b\ell)^{-1/2}\ell^\delta+n^{-\frac{\beta-1}{2(\beta+1)}+\delta}(b\ell)^{(1-\delta)/4}\nonumber\\
&&+\,n^{-1}b^{\frac{2\beta+1}{4(\beta+1)}-\delta}
\ell^{\frac{4\beta+7}{4(\beta+1)}+5\delta}
\Big)+
o_p\Big(
n^{-\frac{\beta-3}{\beta-1}+\delta}(b\ell)^{1/2}+n^{-\frac{3\beta-1}{4(\beta+1)}+\delta}(b\ell)^{1/2}
\nonumber\\
&&+\,n^{-\frac{\beta(2\beta-3)}{(\beta-1)(2\beta+1)}+\delta}b^{\frac{1}{2}}
\ell^{\frac{1}{2}+\frac{2(\beta+3)}{(\beta-1)(2\beta+1)}}+
n^{-\frac{4\beta+5}{4(\beta+1)}+\delta}b^{\frac{1}{2}}
\ell^{\frac{\beta+2}{\beta+1}}+n^{-\frac{2(\beta+1)}{2\beta+1}+\delta}b^{\frac{1}{2}}
\ell^{\frac{4\beta+7}{2(2\beta+1)}}\nonumber\\
&&+\,n^{-\frac{4\beta^2+3\beta+1}{2(2\beta+1)(\beta+1)}+\delta}
b^{\frac{1}{2}}\ell^{\frac{3\beta+4}{2(2\beta+1)}}\Big) \nonumber \\
\label{pf:main3}
\end{eqnarray}
if $\beta\in(5,\infty)$ and $\ell=O(b)$, and
\begin{eqnarray}
\lefteqn{\prob\Big(F^*_n\big(\tilde\xi_n+(b\ell)^{-1/2}x\big)\le p\Big|X_1,\ldots,X_n\Big)}\nonumber\\
&=&\Phi\big(-xf(\xi_p)
/\sigma(\xi_p)\big)+
O_p\Big(\ell^{-1}+\ell^{1/2}n^{-1/2}+(b\ell)^{-1/2}
+n^{-1}b^{\frac{1}{2}-\delta}
\ell^{1+5\delta}\nonumber\\
&&+\,n^{-\frac{1}{2}+\delta}(b\ell)^{(1-\delta)/4}\Big)
+o_p\Big(
n^{-\frac{3}{4}+\delta}(b\ell)^{1/2}+n^{-1+\delta}b^{\frac{1}{2}}\ell\Big)
\label{pf:main4}
\end{eqnarray}
if $\beta=\infty$. Theorem~\ref{thm:main} then follows by (\ref{pf:true}),
(\ref{pf:main3}), (\ref{pf:main4}) and noting that
\begin{eqnarray*}
\prob\Big(F^*_n\big(\tilde\xi_n+(b\ell)^{-1/2}x\big)>p\Big|X_1,\ldots,X_n\Big)
&\le&\prob\Big(\sqrt{b\ell}\big(\xi^*_n-\tilde\xi_n\big)\le x
\Big|X_1,\ldots,X_n\Big)\\
&\le&\prob\Big(F^*_n\big(\tilde\xi_n+(b\ell)^{-1/2}x\big)\ge p\Big|X_1,\ldots,X_n\Big).
\end{eqnarray*}

\section{Conclusion}
\label{sec:conclusion}

In the absence of exact, finite-sample results, accurate estimation of quantiles is essential for implementation of statistical inference procedures. As sample quantiles are non-smooth functionals, conventional bootstrap theory for the smooth function model does not apply to estimation of their distribution. In dependent data settings, with some notable exceptions, little is known about the block bootstrap for distribution estimation of sample quantiles. In this paper we have established a general optimality theory for block bootstrap procedures in such settings under strong mixing conditions, and we have shown that a hybrid block bootstrap is optimal, in the sense of having the fastest convergence rate for distribution estimation. In addition, of course, since the hybrid block bootstrap is based on bootstrap samples of smaller size than the data sample, it provides computational advantage over MBB. How one should choose $(b,\ell)$ in a given application to capture the good theoretical properties of the hybrid block bootstrap requires further consideration. We have provided discussion of an empirical scheme that seems fruitful for this purpose and which will be further developed and refined elsewhere.

\bibliographystyle{plainnat}
\bibliography{peter-ANZJS}

\end{document}